\pdfoutput=1
\RequirePackage{ifpdf}
\ifpdf 
\documentclass[pdftex]{sigma}
\else
\documentclass{sigma}
\fi

\usepackage{mathrsfs}

\renewcommand{\a}{\alpha}
\renewcommand{\b}{\beta}
\newcommand{\ga}{\gamma}
\newcommand{\z}{\zeta}

\newcommand{\la}{\lambda}
\newcommand{\de}{\delta}

\newcommand{\ep}{\varepsilon}

\newcommand{\sigl}{\sigma}

\newcommand{\bfN}{\mathbb{N}}
\newcommand{\bfZ}{\mathbb{Z}}

\newcommand{\Z}{\mathbb{Z}}

\newcommand{\ds}{\displaystyle}
\def\PI{\mbox{\rm P$_{\rm I}$}}
\def\PII{\mbox{\rm P$_{\rm II}$}}
\def\PIII{\mbox{\rm P$_{\rm III}$}}
\def\PIV{\mbox{\rm P$_{\rm IV}$}}
\def\PV{\mbox{\rm P$_{\rm V}$}}
\def\PVI{\mbox{\rm P$_{\rm VI}$}}
\def\PIIIp{$\mbox{P}_{\rm III'}$}

\def\p{Pain\-lev\'e}

\def\bk{B\"ack\-lund}

\def\peqs{\p\ equations}
\def\a{\alpha}
\def\b{\beta}
\def\d{{\rm d}}
\def\H{\mathcal{H}}
\newcommand{\deriv}[3][]{\frac{\d^{#1}{#2}}{{\d{#3}}^{#1}}}
\def\VY{Yablonskii--Vorob'ev polynomials}
\def\Ok{Okamoto polynomials}

\def\L{\mathcal{L}}
\def\eqalign#1{\begin{split}#1\end{split}}
\def\Disc{{\rm Dis}}

\def\fig#1#2{\includegraphics[width=#1]{Figures/#2}}

\newcommand{\Wr}{\operatorname{Wr}}

\numberwithin{equation}{section}

\newtheorem{Theorem}{Theorem}[section]
\newtheorem*{Theorem*}{Theorem}
\newtheorem{Corollary}[Theorem]{Corollary}
\newtheorem{Lemma}[Theorem]{Lemma}

 { \theoremstyle{definition}
\newtheorem{Definition}[Theorem]{Definition}

\newtheorem{Remark}[Theorem]{Remark}
\newtheorem{Remarks}[Theorem]{Remarks} }

\begin{document}
\allowdisplaybreaks

\newcommand{\arXivNumber}{1609.00495}

\renewcommand{\PaperNumber}{080}

\FirstPageHeading

\ShortArticleName{A Constructive Proof for the Umemura Polynomials of the Third Painlev\'e Equation}

\ArticleName{A Constructive Proof for the Umemura Polynomials\\ of the Third Painlev\'e Equation}

\Author{Peter A.~CLARKSON~$^{\rm a}$, Chun-Kong LAW~$^{\rm b}$ and Chia-Hua LIN~$^{\rm b}$}

\AuthorNameForHeading{P.A.~Clarkson, C.-K.~Law and C.-H.~Lin}

\Address{$^{\rm a)}$~School of Mathematics, Statistics and Actuarial Science, University of Kent,\\
\hphantom{$^{\rm a)}$}~Canterbury, CT2 7NF, UK}
\EmailD{\href{mailto:P.A.Clarkson@kent.ac.uk}{P.A.Clarkson@kent.ac.uk}}

\Address{$^{\rm b)}$~Department of Applied Mathematics, National Sun Yat-sen University,\\
\hphantom{$^{\rm b)}$}~Kaohsiung, Taiwan 804, Taiwan}
\EmailD{\href{mailto:law@math.nsysu.edu.tw}{law@math.nsysu.edu.tw}, \href{mailto:j891033@yahoo.com.tw}{j891033@yahoo.com.tw}}

\ArticleDates{Received June 29, 2023, in final form October 17, 2023; Published online October 25, 2023}

\Abstract{We are concerned with the Umemura polynomials associated with rational solutions of the third Painlev\'e equation. We extend Taneda's method, which was developed for the Yablonskii--Vorob'ev polynomials associated with the second Painlev\'e equation, to give an algebraic proof that the rational functions generated by the nonlinear recurrence relation which determines the Umemura polynomials are indeed polynomials. Our proof is constructive and gives information about the roots of the Umemura polynomials.}


\Classification{33E17; 34M55; 65Q30}

\begin{flushright}
\textit{Dedicated to the memory of John Bryce McLeod $($1929--2014$)$}
 \end{flushright}

\section{Introduction}

 The third \p\ equation (\PIII) {has the form}
 \begin{equation}\label{eq:PIII}
 \deriv[2]{w}{z} = \frac{1}{w}\left(\deriv{w}{z}\right)^{2} - \frac{1}{z}\deriv{w}{z} + \dfrac{\alpha w^2 +\beta}{z} +\gamma w^3 +\dfrac{\delta}{w},
 \end{equation}
where $'={\rm d}/{\rm d}z$ and $\alpha$, $\beta$, $\gamma$ and $\delta$ are arbitrary parameters.
We discuss the Umemura polynomials associated with rational solutions of \eqref{eq:PIII} in the generic case when $\gamma\delta \neq 0$, so we set $\gamma=1$ and~$\delta = -1$, without loss of generality (by rescaling $w$ and $z$ if necessary), and so consider
 \begin{equation}
 \deriv[2]{w}{z} = \frac{1}{w}\left(\deriv{w}{z}\right)^{2} - \frac{1}{z}\deriv{w}{z} + \dfrac{\alpha w^2 +\beta}{z} +w^3 -\dfrac{1}{w}.\label{eq:PIIIn}
 \end{equation}

The six \peqs\ (\PI--\PVI), were discovered by \p, Gambier and their colleagues whilst studying second order ordinary differential equations of the form
\[
\deriv[2]{w}{z}=F\left(z,w,\deriv{w}{z}\right), 
\]
where $F$ is rational in ${\rm d}w/{\rm d}z$ and $w$ and analytic in $z$. The \peqs\ can be thought of as nonlinear analogues of the classical special functions. Indeed, Iwasaki, Kimura, Shimomura and Yoshida \cite{refIKSY} characterize the six \peqs\ as ``the most important nonlinear ordinary differential equations'' and state that ``many specialists believe that during the twenty-first century the \p\ functions will become new members of the community of special functions''. Subsequently this has happened as the \p\ equations are a chapter in the NIST \textit{Digital Library of Mathematical Functions}\ \cite[Section~32]{refDLMF}.

The general solutions of the \peqs\ are transcendental in the sense that they cannot be expressed in terms of known elementary functions and so require the introduction of a~new transcendental function to describe their solution. However, it is well known that \mbox{\PII--\PVI\ } possess rational solutions and solutions expressed in terms of the classical special functions -- Airy, Bessel, parabolic cylinder, Kummer and hypergeometric functions, respectively -- for special values of the parameters, see, e.g., \cite{refPAC06review,refFA82,refGLS02}
and the references therein. These hierarchies are usually generated from ``seed solutions'' using the associated \bk\ transformations and frequently can be expressed in the form of determinants.

Vorob'ev \cite{refVor} and Yablonskii \cite{refYab59} expressed the rational solutions of \PII\
 \begin{equation}\label{eq:PII}
 \deriv[2]{w}{z} = 2 w^3 + zw +\alpha,
 \end{equation}
 {which arise only when ${\alpha\in\bfZ}$}, in terms of special polynomials, now known as the \textit{\VY}, {that} are defined through the
 recurrence relation (a second-order, bilinear differential-difference equation)
\begin{equation}\label{eq:PT.RS.YV2}
Q_{n+1}{Q_{n-1}}= zQ_{n}^2 -4\bigg[
Q_{n}\deriv[2]{Q_{n}}{z}-\left(\deriv{Q_{n}}{z}\right)^{2}\bigg]\end{equation}
with $Q_{0}(z)=1$ and $Q_{1}(z)=z$.
It is clear from the recurrence relation \eqref{eq:PT.RS.YV2} that the ${Q_{n+1}}$ are rational functions, though it is not obvious that they are polynomials since one is dividing by~$Q_{n-1}$ at every iteration. In fact, it is somewhat remarkable that the $Q_{n}$ are polynomials.
Taneda~\cite{refTaneda00}, see also \cite{refFOU}, used an algebraic method to prove that the functions $Q_{n}$ defined by~\eqref{eq:PT.RS.YV2} are indeed polynomials.

Umemura \cite{refUmemura01,refUmemura20}\footnote{The paper \cite{refDLMF} was written by Umemura in 1996, for the proceedings of the conference ``\textit{Theory of nonlinear special functions:the \p\ transcendents}'', held in Montr\'{e}al which was never published.} derived special polynomials with certain rational and algebraic solutions of \PIII, \PV\ and \PVI, see also \cite{refNoumi20,refNOOU}.
Recently there have been further studies of the special polynomials associated with
rational and algebraic solutions of \PIII\ \cite{refAmdeb,refBM,refBMS,refPACpiii,refKMii,refMCB97,refMurata95,refOkamotoP3,refOKSO,refUW98};
a review of rational and algebraic solutions of \peqs\ is given in \cite{refPACcmft}.
Several of these papers are concerned with the combinatorial structure and determinant representation of the polynomials, often related to the Hamiltonian structure and affine Weyl symmetries of the \peqs. Additionally, the coefficients of these special polynomials have some interesting, indeed somewhat mysterious, combinatorial properties \cite{refUmemura98,refUmemura01,refUmemura20}.

These special polynomials 
 arise in several applications.
For example, the Umemura polynomials associated with rational solutions of \PIII\ and \PV\ arise as multivortex solutions of the complex sine-Gordon equation \cite{refBP98,refBGG03,refBG03,refOB05},
and in MIMO wireless communication systems \cite{refCCBC}.

We emphasize that the fact that the nonlinear recurrence relation \eqref{eq:PT.RS.YV2} generates polynomials also follows from the $\tau$-function theory associated with the theory of \p\ equations. The $\tau$-functions are in general entire functions. It can be shown that for \PII\ with $\a=m$, the associated $\tau$-function is
\[\tau_m(z)=Q_{m}(z) \exp\left(-\frac{z^3}{24}\right).\]
Consequently, the rational function $Q_m(z)$ has to be a polynomial.
Taneda \cite{refTaneda00} and Fukutani, Okamoto and Umemura \cite{refFOU} independently gave a direct algebraic proof, which is one of the first studies of nonlinear recurrence relations for polynomials. In particular, Taneda \cite{refTaneda00} defined a Hirota-like operator
\[\L(f)=f \deriv[2]{f}{z}- \left(\deriv{f}{z}\right)^{2},\] and showed that if $f(z)$ is a polynomial in $z$,
and $g=zf^2-4 \L(f)$, then $f$ divides $2z g^2-4 \L(g)$. Hence if $f(z)=Q_{m-1}(z)$, then $g(z)=Q_m(z) Q_{m-2}(z)$ and
\[
2z g^2-4 \L(g)=Q_m^2 Q_{m-3}Q_{m-1}+Q_{m-2}^2\big[zQ_m^2-4 \L(Q_m)\big],
\]
so that ${Q_{m-1}}$ divides $zQ_m^2-4 \L(Q_m)$, implying that $Q_{m+1}$ is a polynomial. This is based on the assumption that each $Q_m$ has simple zeros (implying that $Q_m$ and $Q_{m-1}$ have no common zeros),
which in turn can be proved using another identity derived from $\PII$,
\[
\deriv{Q_{m+1}}{z} Q_{m-1}-Q_{m+1}\deriv{Q_{m-1}}{z}=(2m+1)Q_m^2,
\]
which is proved in \cite{refFOU,refTaneda00}, see also \cite{refKanOch}.

{In this paper, we are concerned with $\PIII$ \eqref{eq:PIIIn}. In this case the recurrence relation is } 
\begin{equation}\label{eq:Srecrel}
S_{n+1}S_{n-1} =-z\bigg[S_{n}\deriv[2]{S_n}{z}-\left(\deriv{S_n}{z}\right)^{2}\bigg] -S_{n} \deriv{S_n}{z}+(z+\mu)S_{n}^2,
\end{equation}
{where $\mu$ is a complex parameter; see Theorem \ref{thm:Srecrel} below. The objective is to extend Taneda's method to prove directly and constructively that the rational functions $S_n(z;\mu)$ defined by \eqref{eq:Srecrel} are indeed polynomials. }

 Note that in \eqref{eq:Srecrel}, there is one more term $S_n \deriv{S_n}{z}$, and $z$ in the main term implies that the root $z=0$ of $S_{n}(z;\mu)$, if exists, will accumulate. To employ Taneda's {method}, we define another Hirota-like operator
\[
{\L_z(f)=f\deriv[2]{f}{z}-\left(\deriv{f}{z}\right)^{2}+\frac{f}{z}\deriv{f}{z}.}
\]
Also we need {one more} identity. We find that it is suitable to use the fourth order differential equation satisfied by $S_{n}(z;\mu)$ given
in \cite{refPACpiii}.
This fourth order equation comes from 
the second-order, second-degree equation, often called the
\textit{\p\ $\sigma$-equation}, or \textit{Jimbo--Miwa--Okamoto equation}, satisfied by the Hamiltonian associated with \PIII\ given by \cite{refJM,refOkamotoP3}
\begin{align}
\left(z\deriv[2]{\H_n}{z}-\deriv{\H_n}{z}\right)^{2}&+\bigg\{4\left(\deriv{\H_n}{z}\right)^{2}-z^2\bigg\}\left(z\deriv{\H_n}{z}-2\H_n\right)
\nonumber\\ &+4z\bigg[\mu^2-\left(n-\frac12\right)^2\bigg] \deriv{\H_n}{z}
-2z^2\bigg[\mu^2+\left(n-\frac12\right)^2\bigg]=0.\label{eq:S3}
\end{align}
 {Multiplying \eqref{eq:S3} by $1/z^2$ and differentiating with respect to $z$} gives
\[
z^2\deriv[3]{\H_n}{z}-z\deriv[2]{\H_n}{z}+6z \left(\deriv{\H_n}{z}\right)^{2}+(1-8\H_n)\deriv{\H_n}{z}-\frac{1}{2}z^3{+2z\bigg[\mu^2-\left(n+\frac12\right)^2\bigg]}=0,
\]
then letting
\[\H_n(z;\mu)=z\deriv{}{z}\ln S_{n}(z;\mu)-\frac14z^2-\mu z+\frac1{8},\]
gives \cite[p.~9519]{refPACpiii}
 \begin{align*}
 z^2\bigg[S_n\deriv[4]{S_n}{z}-4 \deriv{S_n}{z} \deriv[3]{S_n}{z}+3 \left(\deriv[2]{S_n}{z}\right)^{2}\bigg]&+2z \left(S_n \deriv[3]{S_n}{z}-\deriv{S_n}{z} \deriv[2]{S_n}{z}\right)\nonumber\\ -4z(z+\mu)\bigg[S_n \deriv[2]{S_n}{z}-\left(\deriv{S_n}{z}\right)^{2}\bigg] 
 &-2 S_n \deriv[2]{S_n}{z}+4\mu S_n \deriv{S_n}{z}= 2n (n+1) S_n^2. \label{eq3.5a}
 \end{align*} 
This equation is also instrumental in the analysis of the case when $z=0$ is a root of $S_n(z;\mu)$, see Section~\ref{sec4} below.

Finally, we remark that this is {not} the first paper on the direct proof for Umemura polynomials. In 1999, Kajiwara and Masuda \cite{refKMii} were able to express $S_n(z;\mu)$ in terms of some Hankel determinant of a $n\times n$ matrix of polynomials ({also known as} Schur functions) that can be obtained from an elementary generating function. However, our proof is constructive, giving more information about the order of roots of $S_n(z;\mu)$.
This information was utilized by Bothner, Miller and Sheng \cite{refBM,refBMS} in their study of the asymptotics of the (scaled) poles and roots of the rational solutions in their so-called ``eye-problem''.

In Section~\ref{sec2}, we describe rational solutions of equation \eqref{eq:PIIIn}. In Section~\ref{sec3}, we extend Taneda's algebraic proof for equation \eqref{eq:PT.RS.YV2} to equation \eqref{eq:Srecrel}. {In Section~\ref{sec4}, we discuss $S_n(0;\mu)$ since $z=0$ is the only location where $S_n(z;\mu)$ can have a multiple root}, and in Section~\ref{sec5}, we discuss our results.

\section[Rational solutions of P\_III]{Rational solutions of $\boldsymbol{{\rm P}_{\rm III}}$}\label{sec2}
The {classification} of rational solutions of equation \eqref{eq:PIIIn}, which is \PIII\ with $\ga=1$ and $\de=-1$, are given in the following theorem.
 \begin{Theorem} \label{th2.1}
Equation \eqref{eq:PIIIn} has a rational solution if and only if 
$\a+\ep\b= 4n$
 with $n\in\Z$ and $\ep =\pm1$. \end{Theorem}
 \begin{proof} See Gromak, Laine and Shimomura \cite[p.\ 174]{refGLS02}; also \cite{refMCB97,refMurata95}.
 \end{proof}

Umemura \cite{refUmemura01,refUmemura20} derived special polynomials associated with rational solutions of \PIII\ \eqref{eq:PIIIn}, which are defined in Theorem \ref{thm:Umpolysp3}, and states that these polynomials are the analogues of the \VY\ associated with rational solutions of \PII\ \cite{refVor,refYab59} and the \Ok\ associated with rational solutions of \PIV\ \cite{refOkamotoP2P4}.

\begin{Theorem}\label{thm:Umpolysp3}
{Suppose that $T_{n}(z;\mu)$ satisfies the {recurrence} relation
\begin{equation}\label{eq:Umemrecrel} z T_{n+1}T_{n-1}
=-z\bigg[T_{n}\deriv[2]{T_{n}}{z}-\left(\deriv{T_{n}}{z}\right)^{2}
\bigg] -T_{n}\deriv{T_n}{z} +(z+\mu)T_{n}^2\end{equation} with
$T_{-1}(z;\mu)=1$ and $T_0(z;\mu)=1$. Then
\begin{gather*}
w_{n}(z;\mu)\equiv w(z;\a_{n},\b_{n})
=\frac{T_{n}(z;\mu-1) T_{n-1}(z;\mu)}{T_{n}(z;\mu) T_{n-1}(z;\mu-1)}
\equiv 1+ \deriv{}{z}\log \frac{{T_{n-1}}(z;\mu-1)}{z^n T_{n}(z;\mu)}
\end{gather*}
satisfies \PIII\
\eqref{eq:PIIIn} with $\a_{n}=2n+2\mu-1$ and
$\b_{n}=2n-2\mu+1$.}
\end{Theorem}

\begin{proof} See Umemura \cite{refUmemura01,refUmemura20}; also \cite{refPACpiii,refKMii}. \end{proof}

{We note that $T_{n}(z;\mu)$ are polynomials in $\xi=1/z$}. It is straightforward to determine a~recurrence relation which generates functions $S_{n}(z;\mu)$ which are polynomials in $z$. These are given in the following theorem.
\begin{Theorem}{\label{thm:Srecrel}
Suppose that $S_{n}(z;\mu)$ satisfies the {recurrence} relation \eqref{eq:Srecrel}, i.e., 
\begin{equation*}
S_{n+1}S_{n-1} =-z\bigg[S_{n}\deriv[2]{S_n}{z}-\left(\deriv{S_n}{z}\right)^{2}\bigg] -S_{n} \deriv{S_n}{z}+(z+\mu)S_{n}^2
\end{equation*}
with $S_{-1}(z;\mu)=S_0(z;\mu)=1$. Then
\[
\eqalign{ w_{n}=w(z;\a_{n},\b_{n})
&=\frac{S_{n}(z;\mu-1) S_{n-1}(z;\mu)}{S_{n}(z;\mu) S_{n-1}(z;\mu-1)}
\equiv 1+\deriv{}{z}\log \frac{S_{n-1}(z;\mu-1)}{S_{n}(z;\mu)}
}\]
satisfies \PIII\ \eqref{eq:PIIIn} with $\a_{n}=2n+2\mu-1$ and
$\b_{n}=2n-2\mu+1$.}\end{Theorem}

\begin{proof} See Clarkson \cite{refPACpiii} and Kajiwara \cite{refKajiwara}; see also Kajiwara and Masuda \cite{refKMii}.
\end{proof}

\begin{Remarks}\quad
\begin{enumerate}\itemsep=0pt
\item[(1)] The rational solutions of \PIII\
\eqref{eq:PIIIn} lie on the lines $\a+\ep\b=4n$, with $\ep=\pm1$, in the $\a$-$\b$ plane.
For any $n \in \mathbb{N}\cup\{0\}$, if $\a_n = 2n +2\mu -1$\ and\ $\b_n = 2n-2\mu+1$ with $\mu\in\mathbb{C}$, then~$\a_n +\b_n =4n$.
\item[(2)] The polynomials $S_{n}(z;\mu)$ and $T_{n}(z;\mu)$, defined by \eqref{eq:Srecrel} and \eqref{eq:Umemrecrel}, respectively, are
related through $S_{n}(z;\mu)=z^{n(n+1)/2}T_{n}(z;\mu)$.
Further $S_{n}(z;\mu)$, {also called Umemura polynomials (for \PIII)}, have the symmetry property $S_{n}(z;\mu)=S_{n}(-z;-\mu)$.
\item[(3)] It is trivial to see that each Umemura polynomial $S_n(z;\mu)$ is monic, and $\deg S_n=\frac12 n(n+1)$ for $n\in \bfN$.

\item[(4)] The Umemura polynomials $S_{n}(z;\mu)$ also arise in the description of algebraic solutions of the special case of \PV\ when $\ga\not=0$ and $\de=0$, i.e.,
\[
\deriv[2]{u}{\z} = \left(\frac{1}{2u} + \frac{1}{u-1}\right)\left(\deriv{u}{\z}\right)^{2} -
\frac{1}{\z}\deriv{u}{\z} + \frac{(u-1)^2}{\z^2}\left(\a u+\frac{\b}{u}\right) +
\frac{\ga u}{\z},\]
when
\begin{gather*} (\a,\b,\ga)=\bigg(\frac12\mu^2,-\frac12\left(n-\frac12\right)^2,-1\bigg),\quad\text{or}\quad
(\a,\b,\ga)=\bigg(\frac12\left(n-\frac12\right)^2,-\frac12\mu^2,1\bigg),
\end{gather*}
see \cite{refPACpv,refPACcmft,refPAC23}, which is known to be equivalent to \PIII\ \eqref{eq:PIIIn}, cf.~\cite[Section~34]{refGLS02}.

\item[(5)] Letting $w(z)=u(\z)/\sqrt{\z}$ with
$\z=\frac14z^2$, in \PIII\ \eqref{eq:PIIIn} yields
\begin{equation}\label{eq:PT.INT.aPIII}
\deriv[2]{u}{\z} = \frac{1}{u}\!\left(\deriv{u}{\z}\right)^{2} -
\frac{1}{\z}\deriv{u}{\z} +\frac{\a u^2}{2\z^2} +
\frac{\b}{2\z} +\frac{u^3}{\z^2} - \frac{1}{u},
\end{equation}
which is known as \PIIIp\ (cf.\ Okamoto \cite{refOkamotoP3}) and is frequently used to determine properties of solutions of \PIII.
However, \PIIIp\ \eqref{eq:PT.INT.aPIII} has algebraic solutions rather than rational solutions~{\cite{refBM22,refMCB97,refMurata95}}.
\end{enumerate}\end{Remarks}

Kajiwara and Masuda \cite{refKMii} derived representations of rational
solutions for \PIII\ \eqref{eq:PIIIn} in the form of determinants,
which are described in the following theorem.

\begin{Theorem}{Let $p_k(z;\mu)$ be the polynomial defined by
\begin{equation*}
\sum_{j=0}^\infty p_j(z;\mu)\la^j =(1+\la)^{\mu}\exp(z\la)
\end{equation*}
with $p_j(z;\mu)=0$ for $j<0$, and $\tau_{n}(z)$ for $n\geq1$, be the
$n\times n$ determinant
\begin{equation*}
\tau_{n}(z;\mu)=\Wr (p_1(z;\mu), p_{3}(z;\mu),
\dots, p_{2n-1}(z;\mu) ),\end{equation*}
where $\Wr(\phi_1,\phi_2,\dots,\phi_n)$ is the Wronskian. Then
\begin{equation*}
w_{n}=w(z;\a_{n},\b_{n},1,-1)=
1+\deriv{}{z}\ln \frac{\tau_{n-1}(z;\mu-1)}{\tau_{n}(z;\mu)}
\end{equation*}
for $n\geq1$, satisfies \PIII\ \eqref{eq:PIIIn} with
$\a_{n}=2n+2\mu-1$ and $\b_{n}=2n-2\mu+1$.}\end{Theorem}

\begin{proof}See Kajiwara and Masuda \cite{refKMii}. \end{proof}

\begin{Remarks}\quad
\begin{enumerate}\itemsep=0pt
\item[(1)] We note that $p_k(z;\mu)=L_k^{(\mu-k)}(-z)$, where $L_k^{(m)}(\z)$
is the \textit{associated Laguerre polynomial}, cf.~\cite[Section~18]{refDLMF}.
\item[(2)] The relationship between the polynomial $S_n (z, \mu)$ and the Wronskian $\tau_{n}(z;\mu)$ is
\[S_n (z; \mu)=c_n\tau_{n}(z;\mu),\qquad 
c_n=\prod_{j=1}^n(2j+1)^{n-j}.\]
\item[(3)] In the special case when $\mu=0$, then
\begin{equation} S_n(z;0)=z^{n(n+1)/2},\label{Sn0}\end{equation}
which is straightforward to show by applying induction to \eqref{eq:Srecrel} with $\mu=0$.
\item[(4)] In the special case when $\mu=1$, then 
\begin{equation*}
S_n (z; 1) = z^{n(n-1)/2} \theta_n(z),\end{equation*}
where $\theta_n(z)$ is the \textit{Bessel polynomial}, sometimes known as the \textit{reverse Bessel polynomial}, given by
\begin{equation*}\theta_n(z)=\sqrt{\frac{2}{\pi}} z^{n+1/2} e^z K_{n+1/2}(z)\equiv\frac{n!}{(-2)^n}L_n^{(-2n-1)}(2z)\end{equation*}
with $K_{\nu}(z)$ the \textit{modified Bessel function}, cf.\ {\cite{refBurch,refCarl,refGrosswald,refKF}}, which arise in the description of point vortex equilibria \cite{refON16}. {We note that Bessel functions also arise in the description of special function solutions of \PIII, see Theorem \ref{thm:P3sf}.}
\end{enumerate}
\end{Remarks}

The recurrence relation \eqref{eq:Srecrel} is nonlinear, so in general, there is no guarantee that the rational function $S_{n+1}(z;\mu)$ thus derived is a polynomial (since one is dividing by $S_{n-1}(z;\mu)$), as was the case for the recurrence relation \eqref{eq:PT.RS.YV2}. However, the \p\ theory guarantees that this is the case through an analysis of the $\tau$-function. A few of these Umemura polynomials $S_n(z;\mu)$, with $\mu$ an arbitrary complex parameter, are given in Table~\ref{tbl21}.
\begin{table}[th] \centering
$
\begin{aligned}
 S_1(z;\mu)={}& z+\mu,\\
 S_2(z;\mu) ={}& \xi^3-\mu,\\
 S_3(z;\mu) ={}& \xi^6-5\mu\xi^3+9\mu\xi-5\mu^2,\\
 S_4(z;\mu) ={}& \xi^{10} -15\mu\xi^7+63\mu\xi^5-225\mu\xi^3+315\mu^2\xi^2 
-175\mu^3\xi+36\mu^2,\\
 S_5(z;\mu) ={}& \xi^{15}-35\mu\xi^{12}+252\mu\xi^{10}+175{\mu}^{2}\xi^{9}-2025\mu\xi^{8} 
+945{\mu}^{2}\xi^{7}\\
&{} -1225\mu (\mu^2-9)\xi^{6}-26082{\mu}^{2}\xi^{5} 
+33075{\mu}^{3}\xi^{4}-350{\mu}^{2} \big(35{\mu}^{2}+36\big) \xi^{3}
\\&{} +11340{\mu}^{3}\xi^{2} -225{\mu}^{2} \big(49\mu^2-36\big) \xi+7{\mu}^{3} \big(875{\mu}^{2}-828\big).
\end{aligned}
$
 \caption{The first few Umemura polynomials $S_n(z;\mu)$, with $\xi=z+\mu$.}\label{tbl21}
 \end{table}

It is straightforward to determine when the roots of $S_n(z;\mu)$ coalesce using discriminants of polynomials.
\begin{Definition}Suppose that \[f(z)=z^m+a_{m-1}z^{m-1}+\dots+a_1z+a_0,\] is a monic polynomial of degree $m$ with roots $\a_1,\a_2,\dots,\a_m$, so
\[
f(z)=\prod_{j=1}^m(z-\a_j).
\]
 Then the \textit{discriminant} of $f(z)$ is
\begin{equation*}\Disc(f) = 
\prod_{1\leq j<k\leq m}(\a_j-\a_k)^2.\end{equation*}
Hence the polynomial $f$ has a multiple root when $\Disc(f)=0$. 
\end{Definition}

\begin{table}[t]
\centering
$
\begin{aligned}
\Disc(S_2(z;\mu))&=-3^3\mu^2,\\
\Disc(S_3(z;\mu))&=3^{12}5^5\mu^6\big(\mu^2-1\big)^2,\\
\Disc(S_4(z;\mu))&=3^{27}5^{20}7^7\mu^{14}\big(\mu^2-1\big)^6\big(\mu^2-4\big)^2,\\
\Disc(S_5(z;\mu))&=3^{66}5^{45}7^{28}\mu^{26}\big(\mu^2-1\big)^{14}
\big(\mu^2-4\big)^6\big(\mu^2-9\big)^2,\\
\Disc(S_6(z;\mu))&=-3^{147}5^{80}7^{63}11^{11} \mu^{44}\big(\mu^2-1\big)^{26}\big(\mu^2-4\big)^{14}\big(\mu^2-9\big)^6 \big(\mu^2-16\big)^2.
\end{aligned}
$
\caption{The discriminants of the Umemura polynomials $S_n(z;\mu)$.}\label{tbl22}
\end{table}

The discriminants of the first few Umemura polynomials $S_n(z;\mu)$ are given in Table \ref{tbl22}. From this we see that $S_2(z;\mu)$ has multiple roots when $\mu=0$, $S_3(z;\mu)$ has multiple roots when~${\mu=0,\pm1}$,~$S_4(z;\mu)$ has multiple roots when $\mu=0,\pm1,\pm2$, $S_5(z;\mu)$ has multiple roots when~${\mu=0,\pm1,\pm2,\pm3}$, and $S_6(z;\mu)$ has multiple roots when $\mu=0,\pm1,\pm2,\pm3,\pm4$. Further the multiple roots occur at $z=0$. This leads to the following theorem.

\begin{Theorem}
The discriminant of the polynomial $S_{n}(z;\mu)$ is given by
\[|\Disc(S_n)|=\prod_{{j=0}}^{{n-1}}(2j+1)^{(2j+1)(n-j)^2}\prod_{k=-n}^n(\mu-k)^{c_{n-|k|}}, \]
where $c_n=
\frac16n^3+\frac14n^2-\frac16n
-\frac18[1-(-1)^n]$ and $\Disc(S_n)<0$ if and only $n=2\ \mbox{mod}\ 4$.
Further the polynomial $S_n(z;\mu)$ has multiple roots at $z=0$ when $\mu=0,\pm1,\pm2,\dots,\pm(n-2)$.
\end{Theorem}

\begin{proof} See Amdeberhan \cite{refAmdeb}.\end{proof}

\begin{Theorem}\label{thm:P3sf}{Equation \eqref{eq:PIII} has solutions expressible in terms of Bessel functions if and only~if
$\a+\ep\b= 4m-2$
with $m\in\Z$ and $\ep =\pm1$.}\end{Theorem}

\begin{proof}See Gromak, Laine and Shimomura \cite[Section~35]{refGLS02}; also \cite{refMW,refUW98}.
\end{proof}

Plots of the roots of the polynomials $S_n(z;\mu)$ for various $\mu$ are given in \cite{refPACpiii}.
Initially for~$\mu$ sufficiently large and negative, the $\frac12n(n+1)$ roots of $S_n(z;\mu)$ form an approximate triangle with $n$ roots on each side. Then as $\mu$ increases, the roots in turn coalesce and eventually for~$\mu$ sufficiently large and positive they form another approximate triangle, similar to the original triangle, though with its orientation reversed. {As shown in Theorem \ref{Thm49} below, as $|\mu|\to\infty$ the roots of $S_n(z;\mu)$ tend to ``triangular structure'' of the roots of the Yablonskii--Vorob'ev polynomial $Q_n(z)$ which arise in the description of the rational solutions of \PII\ \eqref{eq:PII}.}

Bothner, Miller and Sheng \cite{refBM,refBMS}
study numerically how the distributions of poles and zeros of the rational solutions of \PIII\ \eqref{eq:PII} behave as $n$ increases and how the patterns vary with $\mu\in\mathbb{C}$ (note that they use a different notation to our notation).

{It is well known that \PII\ \eqref{eq:PII} arises as the coalescence limit of \PIII,
cf.\ \cite{refInce}. If in \PIII\ \eqref{eq:PIIIn}, we let
\[w(z;\a,\b)=1+\ep u(\z;a),\qquad z={\frac{\z}{\ep}+\frac{4}{\ep^3}},\qquad
\a=2a-\frac{8}{\ep^3},\qquad\b=2a+\frac{8}{\ep^3},\]
then $u(\z;a)$ satisfies
\[\deriv[2]{u}{\z}=2u^3+\z u+a +\ep\bigg\{\left(\deriv{u}{\z}\right)^{2} -u^4 +\frac12\z u^2+au\bigg\}+\mathcal{O}\big(\ep^2\big).\]
Hence in the limit as $\ep\to0$,
\eqref{eq:PIIIn} coalescences to \PII\ \eqref{eq:PII}. In the following theorem, it is shown that the Yablonskii--Vorob'ev polynomial $Q_n(\z)$ arises as
the coalescence limit of the polynomial~$S_n(z;\mu)$ in an analogous way, see also \cite{refPACcmft,refCH08}.

\begin{Theorem}\label{Thm49}
The Yablonskii--Vorob'ev polynomial $Q_n(\z)$ arises as the coalescence limit of the polynomial $S_n(z;\mu)$ given by
\[
Q_{n}(\z)= \lim_{\ep\to0}\left\{\ep^{n(n+1)/2}S_n\left({\frac{\z}{\ep}+\frac{4}{\ep^3}};-\frac{4}{\ep^3}\right)\right\}.
\]
\end{Theorem}

\begin{proof}Since $S_n(z;\mu)$ {satisfies the recurrence relation} \eqref{eq:Srecrel}, then {making the transformation}
\begin{equation}\label{def:Rn}R_{n}(\z;\ep)=\ep^{n(n+1)/2}S_n\left({\frac{\z}{\ep}+\frac{4}{\ep^3}};-\frac{4}{\ep^3}\right),\end{equation}
{to \eqref{eq:Srecrel} yields the recurrence relation}
\begin{gather*}
R_{n+1}R_{n-1}=-
4\bigg[R_{n}\deriv[2]{R_{n}}{\z}-\left(\deriv{R_{n}}{\z}\right)^{2}\bigg]\! +\z R_{n}^2
-\ep^2\bigg\{\!\z \bigg[R_{n}\deriv[2]{R_{n}}{\z}-\left(\deriv{R_{n}}{\z}\right)^{2}\bigg]\!+R_{n}\deriv{R_{n}}{\z}\!\bigg\}.
\end{gather*}
Hence in the limit as $\ep\to0$,
then this coalescences to the equation
\[
R_{n+1}R_{n-1}=-4\bigg[R_{n}\deriv[2]{R_{n}}{\z}-\left(\deriv{R_{n}}{\z}\right)^{2}\bigg] +\z R_{n}^2,\]
which is the recurrence relation for the Yablonskii--Vorob'ev polynomial
$Q_n(\z)$, recall \eqref{eq:PT.RS.YV2}. Further, since
$S_0(z;\mu)=1$ and $S_1(z;\mu)=z-\mu$ we have $R_{0}(\z)=1=Q_0(\z)$
and ${R_{1}(\z)=\z=Q_1(\z)}$. Thus $Q_n(\z)=R_{n}(\z;0)$, for all $n$, as
required. \end{proof}

 {\begin{Remarks}{\rm
\begin{enumerate}\itemsep=0pt
\item[]
\item[(1)] It is not obvious that $R_n(\z;\ep)$ is a polynomial in $\ep$, as well as a polynomial in $\z$. See Lemma \ref{tha.1} for a proof. We give the first few $R_n$ in Table \ref{tblRn}.
\item[(2)]Masuda \cite[Section~A.2]{refMasuda03} discusses the coalescence limit of 
Umemura polynomials to Yablonskii--Vorob'ev polynomials
through the associated Hamiltonians.
\end{enumerate}}\end{Remarks}}

\begin{table}[t]\centering
$\begin{aligned}
R_{1}(\z;\ep)={}&\z,\\
R_{2}(\z;\ep)={}&\z^3+4,\\
R_{3}(\z;\ep)={}&\z^6+20\z^3-80-36\ep^2\z,\\
R_{4}(\z;\ep)={}&\z \big(\z^{9}+60 \z^{6}+11200\big)-252 \ep^2\z^{2} \big(\z^{3}-20\big)+36\ep^4\big(25 \z^{3}+16\big),\\
R_{5}(\z;\ep)={}&\z^{15}+140 \z^{12}+2800 \z^{9}+78400 \z^{6}-3136000 \z^{3}-6272000
-1008 \ep^2 \z \big(\z^{9}-15 \z^{6}\\ &+2100 \z^{3}+2800\big)
+324 \ep^4\z^{2} \big(25 \z^{6}-1288 \z^{3}-2240\big)-252 \ep^6\big(175 \z^{6}+800 \z^{3}\\&-1472\big)
+129600\ep^8\z.
\end{aligned}
$
 \caption{The first few polynomials $R_n(\z;\ep)$, defined by \eqref{def:Rn}.}\label{tblRn}
\end{table}

\begin{Corollary}{As $|\mu|\to\infty$, the roots of $S_n(z;\mu)$ tend to ``triangular structure'' of the roots of the Yablonskii--Vorob'ev polynomial $Q_n(z)$.}\end{Corollary}}

Using the Hamiltonian formalism for \PIII, it is shown in \cite{refPACpiii} that the polynomials $S_n(z;\mu)$ satisfy {a} fourth order bilinear equation and a sixth order, hexa-linear (homogeneous of degree six) difference equation.

\section{\label{sec3}Application of Taneda's method}
\def\fz{\deriv{f}{z}} \def\fzz{\deriv[2]{f}{z}} \def\fzzz{\deriv[3]{f}{z}}
\def\gz{\deriv{g}{z}} \def\gzz{\deriv[2]{g}{z}} \def\gzzz{\deriv[3]{g}{z}}
\def\hz{\deriv{h}{z}} \def\hzz{\deriv[2]{h}{z}}
 In this section, we use the algebraic method due to Taneda \cite{refTaneda00} to prove that the rational functions $S_n(z;\mu)$ satisfying \eqref{eq:Srecrel} are indeed polynomials, {assuming that all the zeros of $S_n(z;\mu)$ are simple}.

We define an operator $ \L_z$ as follows:
\[
 \L_z (f) = f \fzz - \left(\fz\right)^{2} +\frac{f}{z}\fz.
\]

\begin{Lemma}\label{th2.15}
 Let $f(z)$ and $g(z)$ be arbitrary polynomials. Then
 \begin{enumerate}\itemsep=0pt
 \item[$(a)$]
 $\L_z(k f)=k^2 \L_z(f)$ with $k$ a constant;
 \item[$(b)$] $\L_z(f g)=f^2\L_z(g)+g^2\L_z(f)$;
 \item[$(c)$] If $h = -z\L_z (f) +k(z+\mu) f^2$ with $k$ and $\mu$ constants, then
$f \,|\, z\L_z(h)-2k(z+\mu) h^2$,
where the symbol $|$ means that the right-hand side is divisible by the left-hand side.
 \end{enumerate}
\end{Lemma}
\begin{proof}
 (a) This follows directly from the definition.

(b) We observe that{\samepage
 \begin{align*}
 \L_z (fg)&= fg\deriv[2]{}{z}(fg) - \left[\deriv{}{z}(fg)\right]^2 +\dfrac{fg}{z} \deriv{}{z}(fg)\\
 &= f^2\bigg[g\gzz-\left(\gz\right)^{2} +\frac{g}{z}\gz\bigg] +g^2\bigg[f\fzz-\left(\fz\right)^{2} +\frac{f}{z}\fz\bigg]\\
 &= f^2\L_z (g) +g^2 \L_z (f),
 \end{align*}
so the result is valid.}

 (c) Finally, by definition
 \begin{gather*}
 h = -z\bigg[f \fzz - \left(\fz\right)^{2} +\dfrac{f}{z}\fz\bigg] +k(z+\mu) f^2= z \left(\fz\right)^{2} + f\times \text{(a polynomial)},\\
 \hz = -2f\fzz -z\left(f\fzzz-\fz\fzz\right) +kf^2 +2k(z+\mu)f\fz\\
 \phantom{\hz}= z\fz\fzz + f\times \text{(a polynomial)},\\
 \hzz = -\fz\fzz -3f\fzzz -z\bigg[f\deriv[4]{f}{z}-\left(\fzz\right)^{2}\bigg]\! +4kf\fz+2k(z+\mu)\bigg[ f\fzz+\left(\fz\right)^{2}\!\bigg]\\
 \phantom{\hzz}= z\left(\fzz\right)^{2} -\fz\fzz +2k (z+\mu)\left(\fz\right)^{2} +f\times \text{(a polynomial)}.
 \end{gather*}
 Then we can see
 \begin{align*}
 \L_z (h)&= h \hzz - \left(\hz\right)^{2} +\dfrac{h}{z}\hz\\
 &= {z \left(\fz\right)^{2}\bigg\{z\left(\fzz\right)^{2} -\fz\fzz +2k (z+\mu)\left(\fz\right)^{2}\bigg\} - \left(z\fz\fzz\right)^{2}}\\
 &\quad{}+ z\left(\fz\right)^3 \fzz + f\times \text{(a polynomial)}\\
 &= 2kz(z+\mu)\left(\fz\right)^4 
 +f\times \text{(a polynomial)}.
 \end{align*}
 Since $z\L_z (h) -2k (z+\mu)h^2 = f\times \text{(a polynomial)}$, then
 \begin{equation}
 f \,|\, z\L_z (h) -2k (z+\mu)h^2 \label{eq4.20}
 \end{equation} as required.
\end{proof}

\begin{Theorem}\label{th2.12}
 Suppose $\{ S_n(z;\mu) \}$ is a sequence of rational functions with simple nonzero roots, satisfying \eqref{eq:Srecrel}, with $S_{-1}(z;\mu)=S_0(z;\mu)=1$. For all $N\in \bfN\cup\{ 0\}$,
 if $z=0$ is not a root of any $S_n(z;\mu)$ for $0\leq n\leq N$, then
 \begin{enumerate}\itemsep=0pt
 \item[$(a)$] $S_{N+1}(z;\mu)$ is a polynomial;
 \item[$(b)$] $S_{N+1}(z;\mu)$ and $S_N(z;\mu)$ do not have a common root.
 \end{enumerate}
\end{Theorem}

\begin{proof} {We first prove part (b).} If $S_{N}(z;\mu)$ and $S_{N-1}(z;\mu)$ have the same root $z_0 \neq 0$, then by~\eqref{eq:Srecrel}, $z_0$ is also a root of
\[S_N \deriv[2]{S_N}{z}-\left(\deriv{S_N}{z}\right)^{2},\] and hence also a root of $\deriv{S_N}{z}(z;\mu)$. This implies $z_0$ is (at least) a double root of $S_{N}(z;\mu)$, which contradicts our assumption about $S_{N}(z;\mu)$.

Part (a) will be shown using induction. First, we have $S_{-1}(z;\mu) = S_{0}(z;\mu) = 1$, then $S_1 (z;\mu) = z + \mu$ and $S_2 (z;\mu) = (z + \mu)^3 -\mu$. Clearly, (a) hold for $n = 0, 1, 2$, when $\mu\neq 0$. We next assume that (a) hold for $n=N-2,N-1,N$ with $N\geq 2$. Then we will prove that the statements also hold for $n = N+1$.

 Let $f$ be $S_{N-1}$. Then $n = N-1$ and $h = S_{N}S_{N-2}$ in Lemma~\ref{th2.15}. Then \eqref{eq4.20} becomes
 \[
 S_{N-1}\,|\, z\L_z (S_{N}S_{N-2}) {+2} (z+\mu)(S_{N}S_{N-2})^2.
 \]
 Hence
 \begin{gather*}
z \left[\L_z (S_{N}S_{N-2}) -\dfrac{2(z+\mu)}{z}(S_{N}S_{N-2})^2\right] \\
\qquad = z\big[S_{N-2}^2 \L_z (S_{N}) +S_{N}^2 \L_z (S_{N-2})\big] -2(z+\mu)(S_{N}S_{N-2})^2\\
 \qquad= S_{N-2}^2\big[z\L_z (S_{N}) -(z+\mu)S^2_{N} \big] +S_{N}^2\big[z\L_z (S_{N-2}) -(z+\mu)S^2_{N-2} \big]\\
 \qquad= S_{N-2}^2\big[z\L_z (S_{N}) -(z+\mu)S^2_{N} \big] -S_{N}^2 S_{N-1}S_{N-3}.
 \end{gather*}
 Then by \eqref{eq4.20} and (b) with $n = N-1$, we have
\[
 S_{N-1}\,|\, {-}z\L_z (S_{N}) +(z+\mu)S^2_{N} = -z\bigg[S_N \deriv[2]{S_N}{z} -\left(\deriv{S_N}{z}\right)^{2}\bigg] -S_N \deriv{S_N}{z}+ (z+\mu)S_N^{2}.
\]
 So, according to \eqref{eq:Srecrel}, $S_{N+1}$ is a polynomial by induction.
\end{proof}

\section[Roots of S\_n(z; mu)]{Roots of $\boldsymbol{S_n(z;\mu)}$}\label{sec4}
In this section we {initially} discuss $S_n(0;\mu)$ since $z=0$ is the only location where $S_n(z;\mu)$ can have a multiple root.
 \begin{Theorem}\label{th3.1}
 Let $\phi_n=S_n(0;\mu)$, and
 \begin{gather*}
 \phi_n':= \frac{\partial S_n}{\partial z}(0;\mu),\qquad \phi_n'':= \frac{\partial^2 S_n}{\partial z^2}(0;\mu),\end{gather*}
 etc. Then for all $n\geq 3$,
 \begin{align}
 \phi_{n+1}&=\frac{\phi_n \phi_{n-1}}{\phi_{n-2}} \left(2\mu^2-2n^2+2n-1-\frac{\phi_n \phi_{n-3}}{\phi_{n-1}\phi_{n-2}}\right);\label{eq3.a}
 \\
 \phi'_{n+1} &= -\frac{\phi_n \phi_{n+2}}{\phi_{n+1}} + \mu \phi_{n+1}. \label{eq3.b}
 \end{align}
 \end{Theorem}

 \begin{proof}
Differentiating \eqref{eq:Srecrel} with respect to $z$ gives
 \begin{align}
 \deriv{S_{n+1}}{z}=\frac{1}{S_{n-1}} \bigg\{& S_n^2+2(z+\mu)S_n \deriv{S_n}{z}-2 S_n \deriv[2]{S_n}{z}
\nonumber\\
&{}+z\bigg(\deriv{S_n}{z} \deriv[2]{S_n}{z}-S_n\deriv[3]{S_n}{z}\bigg)-S_{n+1}\deriv{S_{n-1}}{z}\bigg\}. \label{eq3.2}
 \end{align}
 Substitute $z=0$ into \eqref{eq:Srecrel} and \eqref{eq3.2}. We obtain
 \begin{gather}
 \phi_{n+1}=\frac{\phi_n}{\phi_{n-1}}\left(\mu \phi_n-\phi_n'\right),
\label{eq3.3}
\\
 \phi_{n+1}'=\frac{\phi_n}{\phi_{n-1}} \left( \phi_n+2\mu \phi_n'-2\phi_n''-\frac{\phi_{n-1}'\phi_{n+1}}{\phi_n} \right).
 \label{eq3.4}
 \end{gather}
 Now \eqref{eq3.3} implies that \eqref{eq3.b} is valid. Furthermore, in \cite[p.\ 9519]{refPACpiii}, it was shown that
 \begin{align}
 z^2&\bigg[S_n\deriv[4]{S_n}{z}-4 \deriv{S_n}{z} \deriv[3]{S_n}{z}+3 \left(\deriv[2]{S_n}{z}\right)^{2}\bigg]+2z \left(S_n \deriv[3]{S_n}{z}-\deriv{S_n}{z} \deriv[2]{S_n}{z}\right)\nonumber\\
 &\quad -4z(z+\mu)\bigg[S_n \deriv[2]{S_n}{z}-\left(\deriv{S_n}{z}\right)^{2}\bigg] 
 -2 S_n \deriv[2]{S_n}{z}+4\mu S_n \deriv{S_n}{z}= 2n (n+1) S_n^2. \label{eq3.5}
 \end{align}
 This implies, as $\phi_n$ is not identically zero, that
 \begin{equation}
 2\mu \phi_n'-\phi_n''=n(n+1) \phi_n.
 \label{eq3.13}
 \end{equation}
 Hence by \eqref{eq3.b},
\[
 \phi_n'' = 2\mu \phi_n'-n(n+1)\phi_n=\big[2\mu^2-n(n+1)\big]\phi_n-\frac{2\mu \phi_{n-1}\phi_{n+1}}{\phi_n}.
\]
 Now substitute this equation and \eqref{eq3.b} into \eqref{eq3.4} to obtain, after simplification,
\[
 -\frac{\phi_n \phi_{n+2}}{\phi_{n+1}} = \frac{\phi_n^2}{\phi_{n-1}} \big(2n^2+2n{+1}-2\mu^2\big)+\frac{\phi_{n+1}\phi_n \phi_{n-2}}{\phi_{n-1}^2}.
\]
 Therefore, we have
\[
 \phi_{n+2}=\frac{\phi_n \phi_{n+1}}{\phi_{n-1}} \left( 2\mu^2-2n^2-2n{-1}-\frac{\phi_{n+1} \phi_{n-2}}{\phi_n \phi_{n-1}} \right),
\]
and so \eqref{eq3.a} is also valid.
\end{proof}

 \begin{Corollary}\quad\label{th3.2}
 \begin{enumerate}\itemsep=0pt
 \item[$(a)$] For all $n\in \bfN$,
\[
 \phi_n(\mu)=\mu^{\ga_0^n} \prod_{j=1}^{n-1}(\mu^2-j^2)^{\ga_j^n},
\]
\samepage{where for $0\leq j<k$,
\begin{alignat*}{3}
 &\ga_{2j}^n = \left\lceil \frac{n}{2}\right\rceil -j =k-j \quad&&\mbox{if } \quad n=2k\quad\mbox{ or }\quad 2k-1;&\\
 & \ga_{2j+1}^n = \left\lfloor \frac{n}{2}\right\rfloor -j =k-j \quad &&\mbox{if }\quad n=2k\quad\mbox{ or }\quad 2k+1.&
\end{alignat*}
 \item[$(b)$]When $n\geq 3$, $\phi_n'=\phi_{n-1}g_n(\mu)$, where {$g_n$ is a polynomial of degree $n-1$}.}
 \end{enumerate}
\end{Corollary}

 \begin{Remark}
 \label{rmk}
 Part (a) above means that $z=0$ is a root of $S_n(z;\mu)$ if and only if $\mu =0,\allowbreak\pm 1,\pm 2,\dots,\pm(n-1)$,
 {i.e., $|\mu|$ is an integer strictly less than $n$.}
 In particular, the first few $\phi_n(\mu)$ are
 \begin{align*}
 \phi_1 &= \mu,\\
 \phi_2 &= \mu\big(\mu^2-1\big),\\
 \phi_3 &= \mu^2\big(\mu^2-1\big)\big(\mu^2-4\big),\\
 \phi_4 &= \mu^2\big(\mu^2-1\big)^2\big(\mu^2-4\big)\big(\mu^2-9\big),\\
 \phi_5 &= \mu^3\big(\mu^2-1\big)^2\big(\mu^2-4\big)^2\big(\mu^2-9\big)\big(\mu^2-16\big).
 \end{align*}
\end{Remark}

 \begin{proof} {It is trivial to verify by induction hypothesis, with the help of above and \eqref{eq3.a} that,
 \begin{align*}
 &\phi_{2k}= \mu^k(\mu^2-1)^k \prod_{j=1}^{k-1}\big[\mu^2-(2j)^2\big]^{k-j}\big[\mu^2-(2j+1)^2\big]^{k-j};\\
 &\phi_{2k+1}= \mu^{k+1}(\mu^2-1)^k \prod_{j=1}^{k}\big[\mu^2-(2j)^2\big]^{k+1-j}\big[\mu^2-(2j+1)^2\big]^{k-j}.
 \end{align*}
 This proves (a). Also we have
\[
 {\frac{\phi_{2k+1}}{\phi_{2k}}=\mu \prod_{j=1}^k\big(\mu^2-(2j)^2\big),\qquad
 \frac{\phi_{2k}}{\phi_{2k-1}}=\prod_{j=1}^k \big(\mu^2-(2j-1)^2\big).}
 \]
 Hence by \eqref{eq3.b},
 \[
 \phi_{2k}'=\left( \frac{\mu \phi_{2k}}{\phi_{2k-1}}-\frac{\phi_{2k+1}}{\phi_{2k}}\right) \phi_{2k-1}:=\phi_{2k-1}g_{2k}(\mu),
 \]
 where $g_{2k}$ is a polynomial of degree $2k-1$.
 Similarly, by \eqref{eq3.b} again,
 \[
 \phi_{2k+1}'=\left( \frac{\mu \phi_{2k+1}}{\phi_{2k}}-\frac{\phi_{2k+2}}{\phi_{2k+1}}\right) \phi_{2k}:=\phi_{2k}g_{2k+1}(\mu),
\]
 where $g_{2k+1}$ is a polynomial of degree $2k$. Thus the proof of (b) is complete.}
\end{proof}

\begin{Theorem}\label{thm4.4}
 Fix $m\in {\bfN\cup\{0\}}$. Then for the recurrence relation \eqref{eq:Srecrel} with initial polynomials $S_{-1}=S_{0}=1$, we have
 \begin{enumerate}\itemsep=0pt
 \item[$(a)$] all the non-zero roots of rational functions $S_n(z;\pm m)$ are simple, {for all $n\in\bfN$};
 \item[$(b)$] each $S_n(z;\pm m)$ is a polynomial in $z$, for $n=0,1,\dots,m$.
 \end{enumerate}
 \end{Theorem}

 \begin{proof} We shall make use of the identity \eqref{eq3.5} again. Suppose $z_0$ is a nonzero root of $S_n(z;\mu)$. Then from \eqref{eq3.5},
\[
 3{z_0}\left[\deriv[2]{S_n}{z}(z_0)\right]^2=\deriv{S_n}{z}(z_0)\left[4{z_0} \deriv[3]{S_n}{z}(z_0)+2\deriv[2]{S_n}{z}(z_0){-}4(z_0+\mu)\deriv{S_n}{z}(z_0)\right].
\]
 Hence if $z_0$ is a root of $\ds\deriv{S_n}{z}$, then it also has to be a root of $\ds\deriv[2]{S_n}{z}$. That is, if $z_0$ is not a simple root of $S_n(z;\mu)$, then its order $k\geq 3$. Analyzing on the identity \eqref{eq3.5}, the term
\[S_n\deriv[4]{S_n}{z}-4\deriv{S_n}{z}\deriv[3]{S_n}{z}+3 \left(\deriv[2]{S_n}{z}\right)^{2},\]
 has the zero $z_0$ with order {at least} $2k-4$, while the other terms has order {at least $ 2k-3$}. Therefore, let $S_n(z;\mu)=(z-z_0)^k g(z)$, where $g(z)$ is a polynomial and $g(z_0)\neq 0$. Then there exists a polynomial $h(z)$ such that
\[
 S_n\deriv[4]{S_n}{z}-4\deriv{S_n}{z}\deriv[3]{S_n}{z}+3\left(\deriv[2]{S_n}{z}\right)^{2} = (z-z_0)^{2k-4} \left[(z-z_0)h(z)+6k(k-1)g^2(z)\right].
\]
 But {the expression inside the bracket} must have $z_0$ as a root. This gives a contradiction. We see that every {nonzero} $z_0$ is at most a simple root. {This proves (a). Part (b) follows directly from Remark~\ref{rmk} and Theorem~\ref{th2.12}.}
\end{proof}

\begin{figure}[t]
\begin{tabular}{ccc}
\fig{48mm}{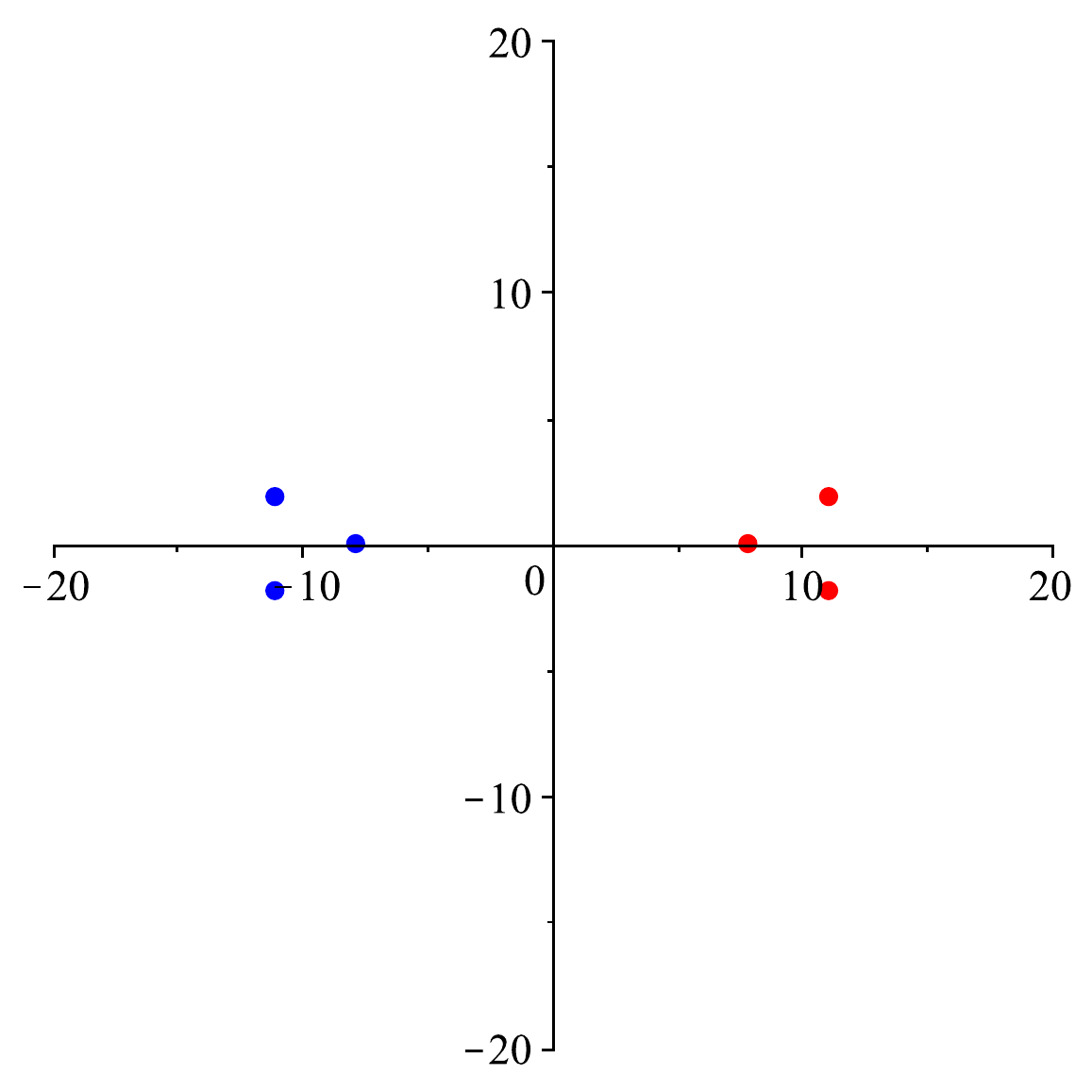} &\fig{48mm}{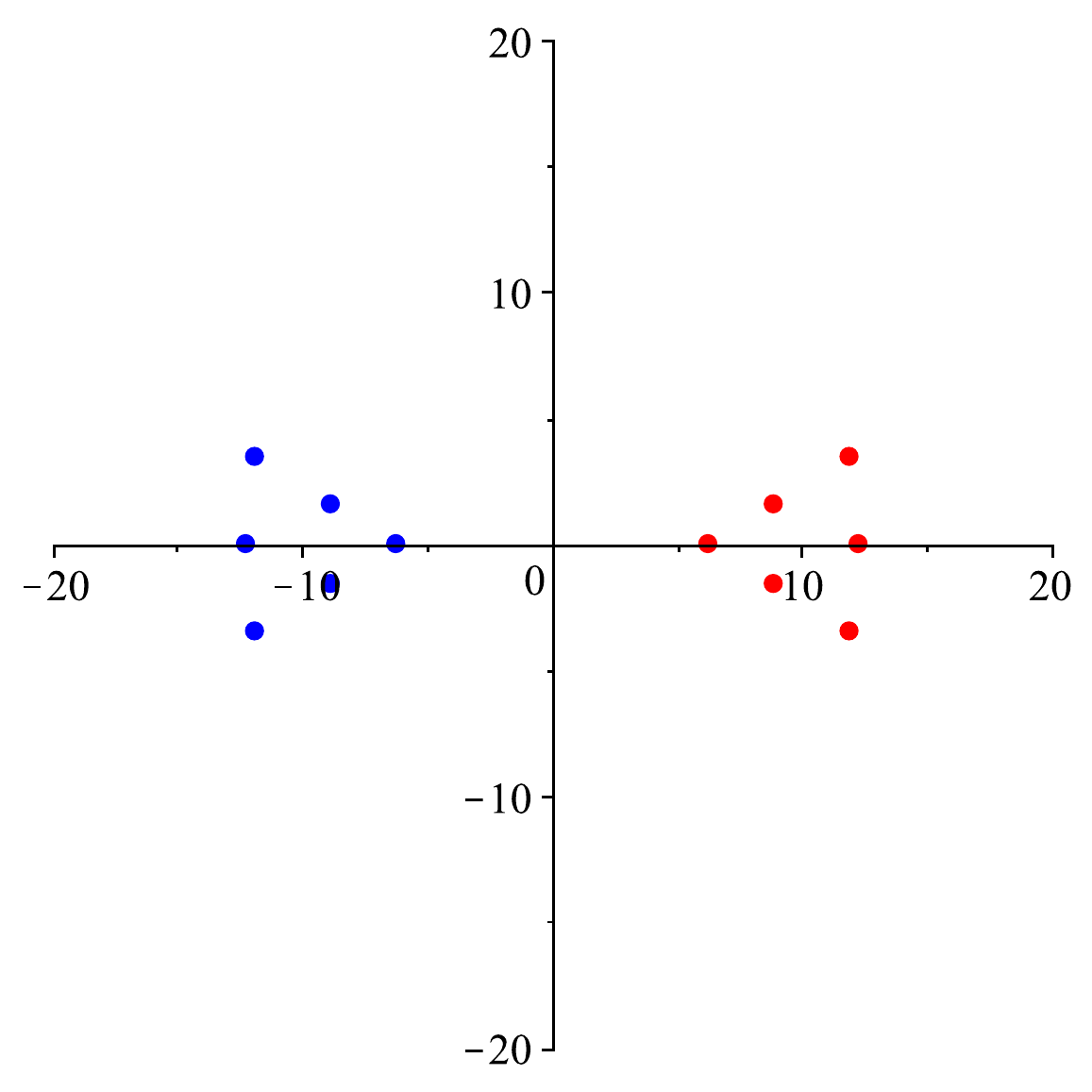} & \fig{48mm}{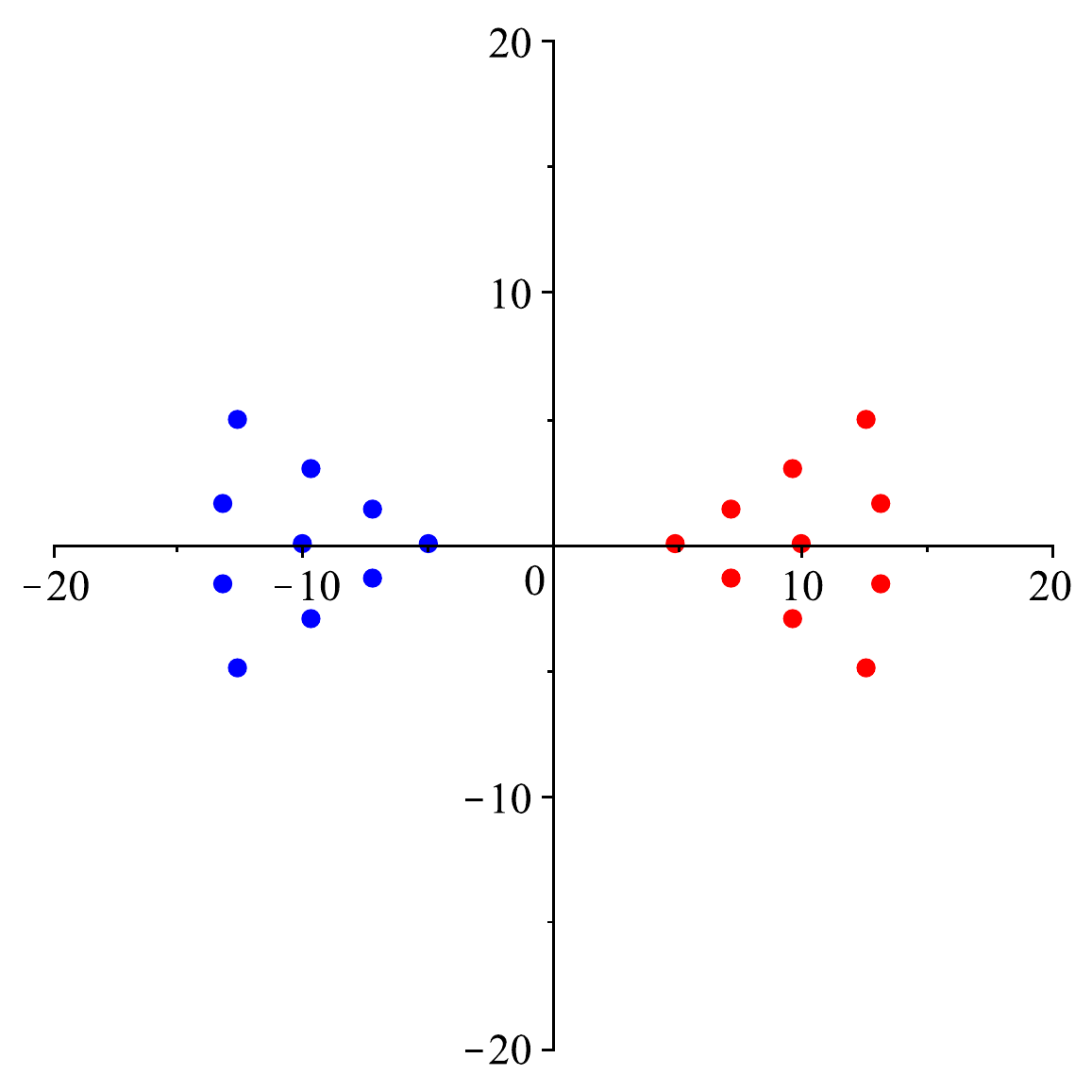}\\
\fig{48mm}{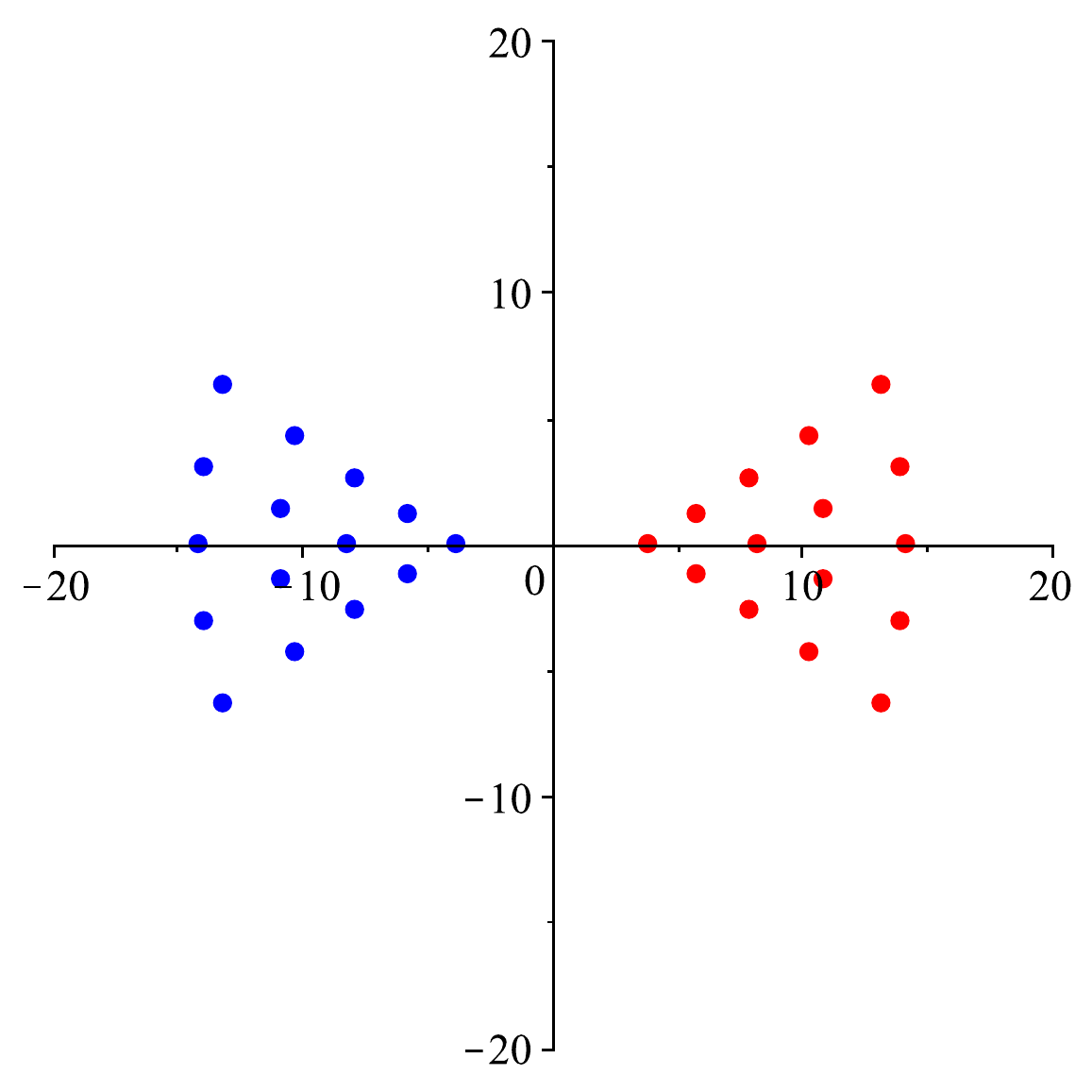} & \fig{48mm}{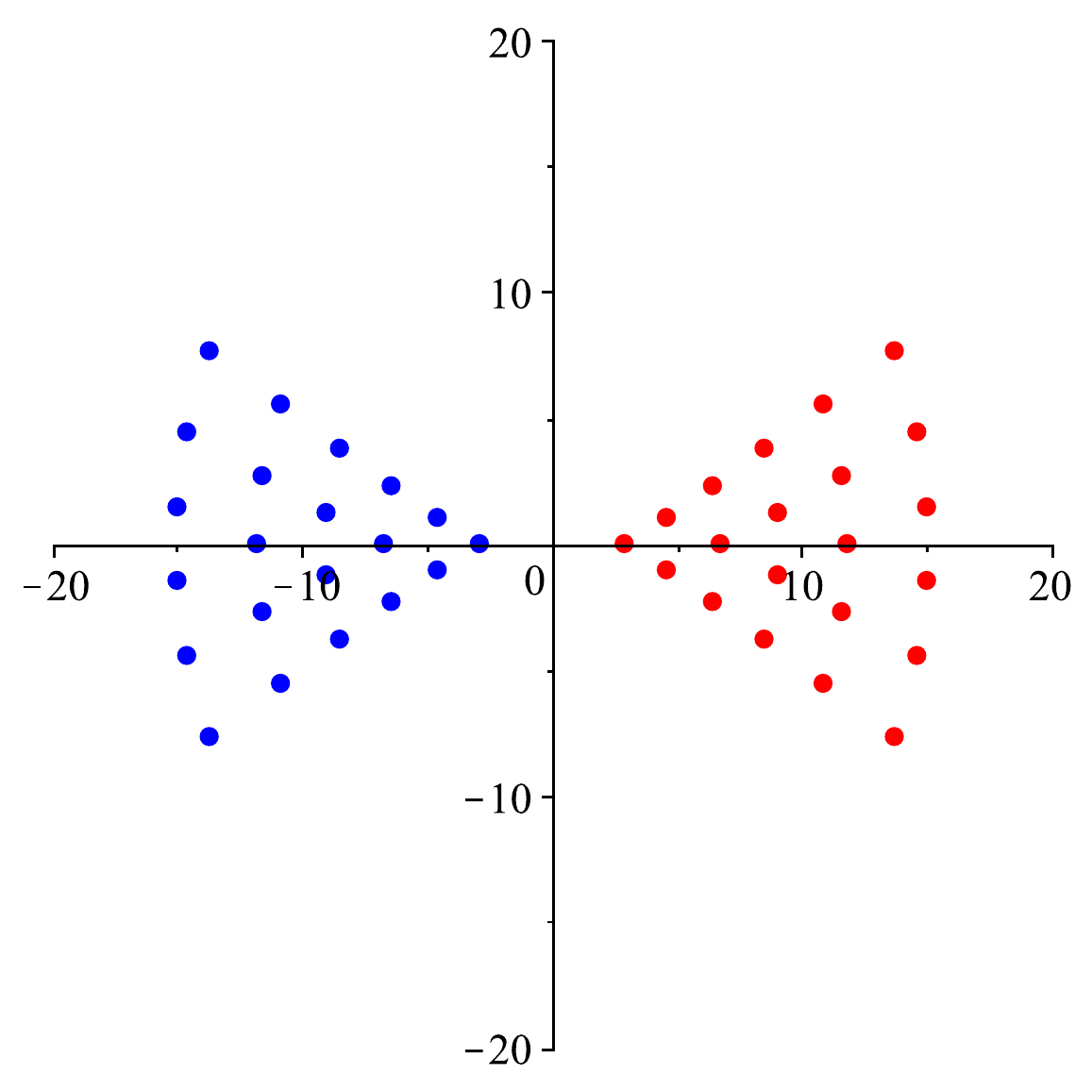} &
\fig{48mm}{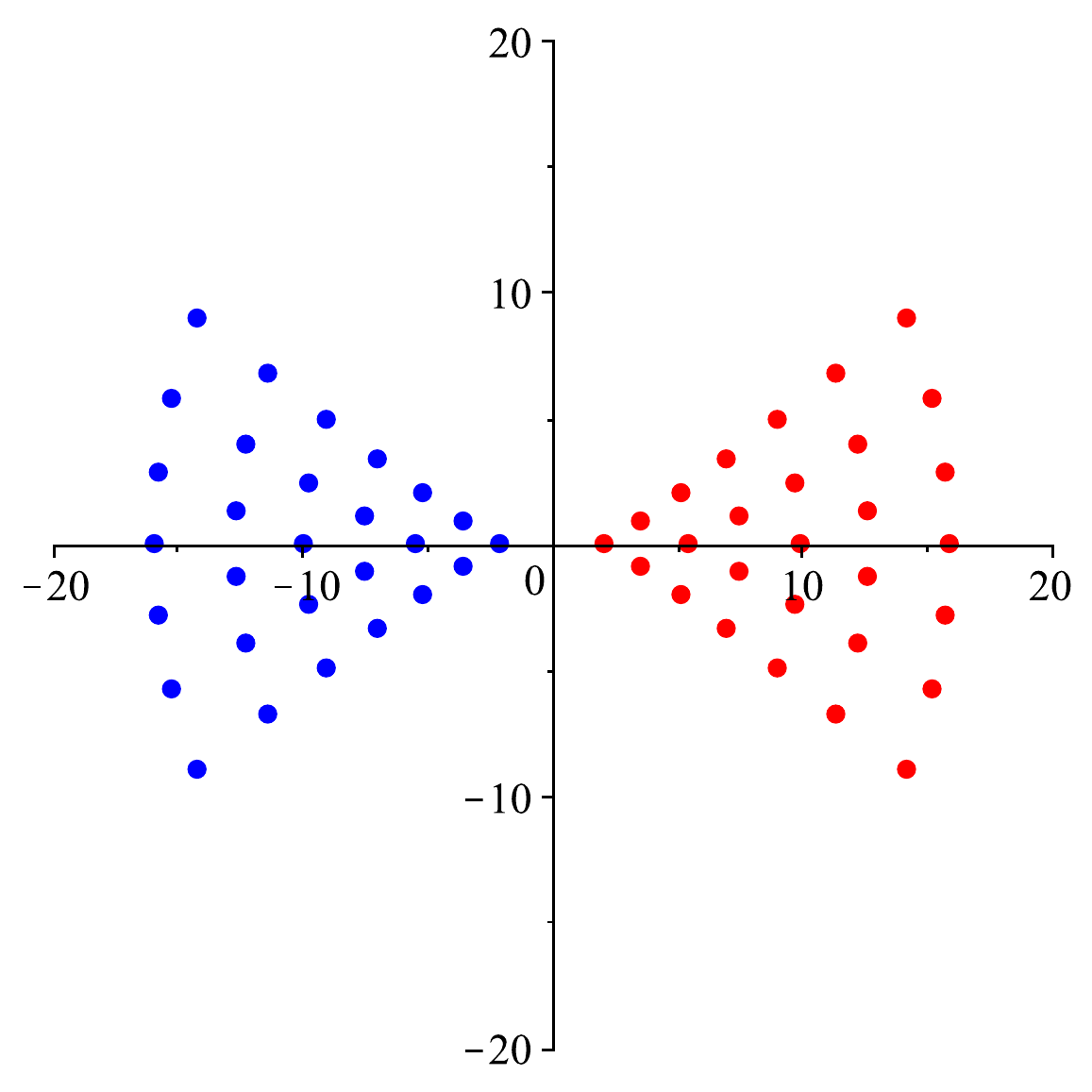}\\
\fig{48mm}{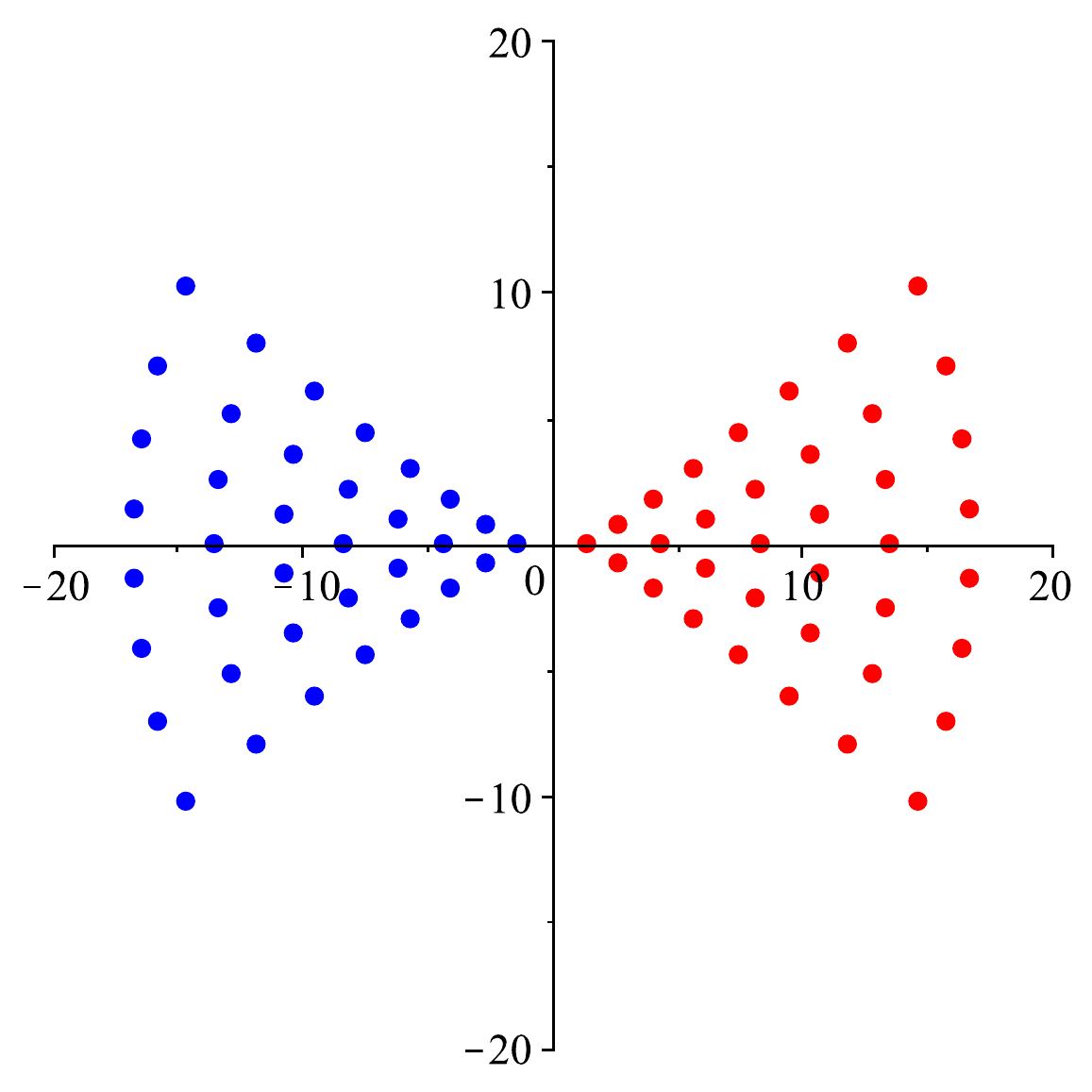} & \fig{48mm}{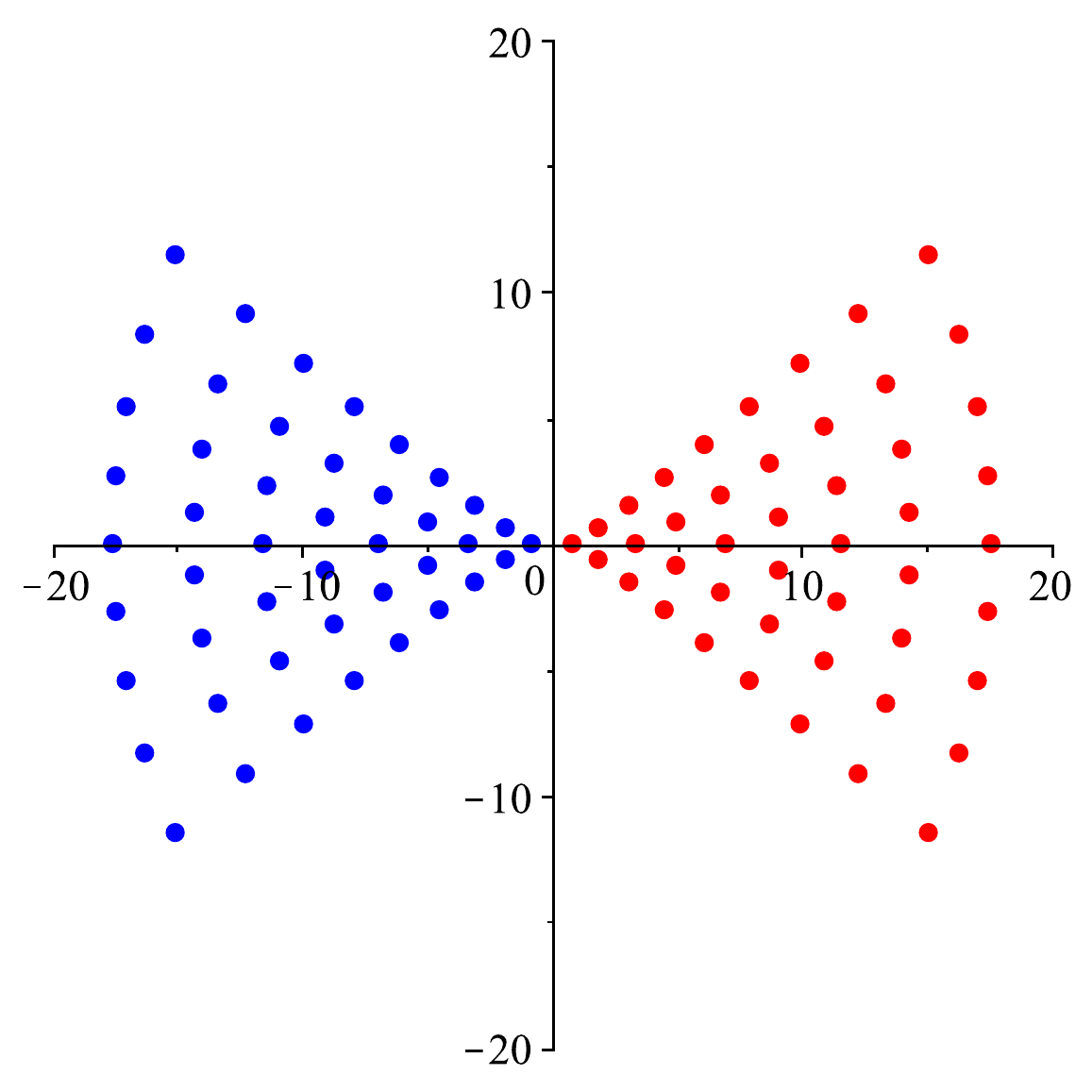} &
\fig{48mm}{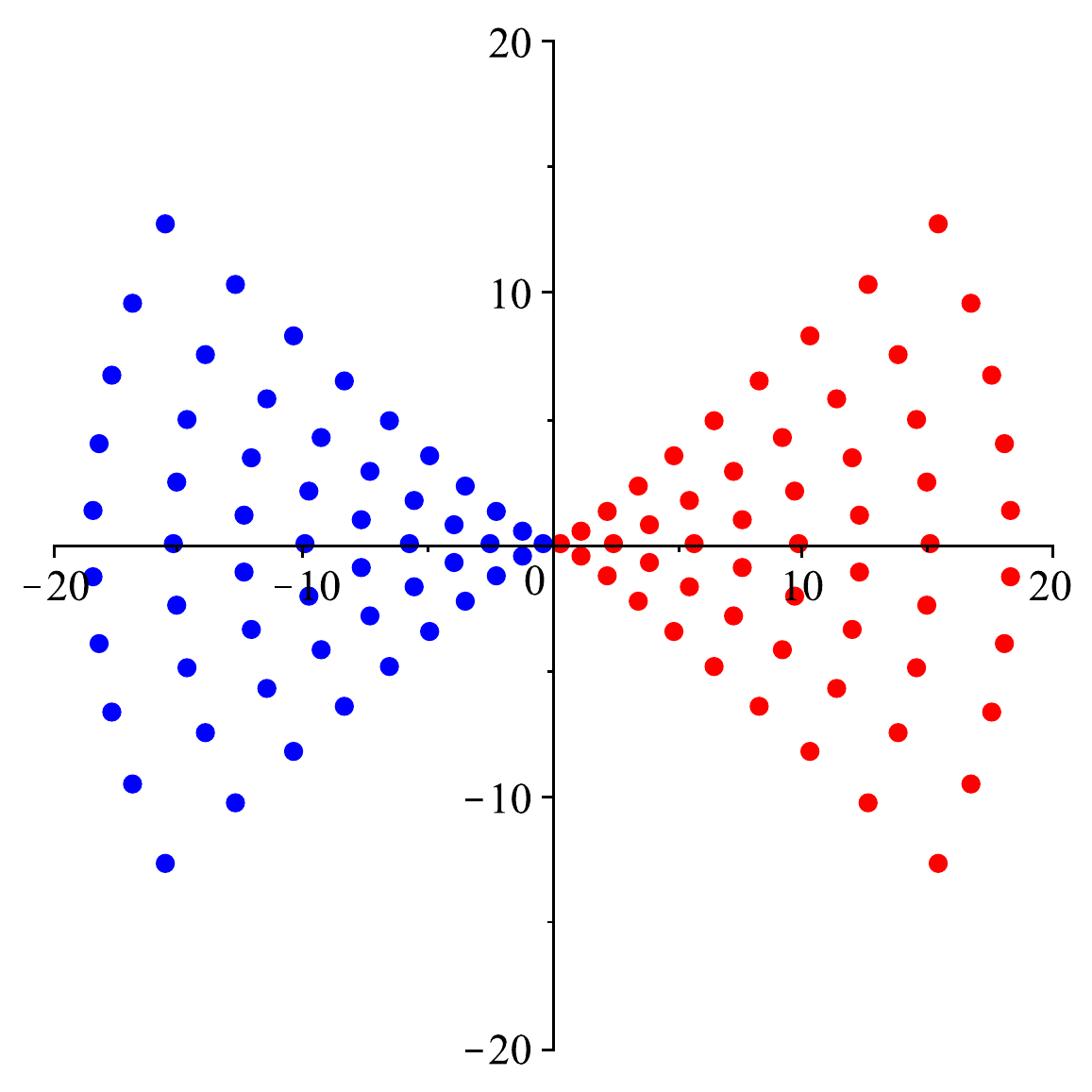}
 \end{tabular}
 \caption{Plots of $S_{n}(z;\mu)$ with $\mu=10$ (blue) and $\mu=-10$ (red), for $n=2,3,\dots,10$. These illustrate the results given in Theorem \ref{thm4.4}. }
 \label{fig:Thm44}
 \end{figure}

These results are illustrated in Figure \ref{fig:Thm44}, where plots of $S_{n}(z;\mu)$ with $\mu=10$ (blue) and $\mu=-10$ (red), for $n=2,3,\dots,10$ are given. Similar figures appear in \cite{refPACpiii}.

 \begin{Lemma} \label{thm4.5}\samepage
 Let $\mu\in \bfZ\setminus \{ 0\}$. Suppose that $S_n(z;\mu)=z^{\sigl} g(z)$, where $g(z)=\sum_{j=0}^k a_j z^j$ {is a~polynomial} $(a_0 \neq 0)$.
 Then
 \begin{enumerate}\itemsep=0pt
 \item[$(a)$] $a_1 =\mu a_0$;
 \item[$(b)$] if {$\sigl=\frac12 \ell (\ell +1)$ with $\ell =n-|\mu|$}, then
$a_2=\frac12\big(\mu^2-\frac{|\mu|}{2\ell +1}\big) a_0$.
 \end{enumerate}
 \end{Lemma}

 \begin{proof} We use the auxilliary identity \eqref{eq3.5} for the proof. First
 \begin{gather}
 \deriv{S_n}{z}=z^{\sigl} \gz+\sigl z^{\sigl-1} g,\nonumber\\
 \deriv[2]{S_n}{z} =  z^{\sigl} \gzz +2\sigl z^{\sigl-1}\gz+\sigl(\sigl-1) z^{\sigl-2} g,\nonumber\\
 \deriv[3]{S_n}{z} =  z^{\sigl} \gzzz +3\sigl z^{\sigl-1}\gzz+3\sigl (\sigl-1)z^{\sigl-2}\gz+\sigl(\sigl-1)(\sigl-2)z^{\sigl-3}g,\nonumber\\
 \deriv[4]{S_n}{z} =  z^{\sigl} \deriv[4]{g}{z} +4\sigl z^{\sigl-1}\gzzz+6\sigl(\sigl-1)z^{\sigl-2}\gzz+4\sigl(\sigl-1)(\sigl-2)z^{\sigl-3}\gz\nonumber\\
 \hphantom{\deriv[4]{S_n}{z} =}{}
 +\sigl(\sigl-1)(\sigl-2)(\sigl-3)z^{\sigl-4}g.\label{SnSn4}
 \end{gather}
 Express \eqref{eq3.5} as
 \begin{gather}
2n(n+1)S_n^2+4z^2\bigg[S_n\deriv[2]{S_n}{z}-\left(\deriv{S_n}{z}\right)^{2}\bigg]\nonumber\\
\qquad = z^2\bigg[S_n \deriv[4]{S_n}{z}-4 \deriv{S_n}{z} \deriv[3]{S_n}{z}+3 \left(\deriv[2]{S_n}{z}\right)^{2}\bigg]+2z \left(S_n \deriv[3]{S_n}{z}-\deriv{S_n}{z} \deriv[2]{S_n}{z}\right)\nonumber\\
 \qquad\quad-4\mu z\bigg[S_n \deriv[2]{S_n}{z}-\left(\deriv{S_n}{z}\right)^{2}\bigg] 
 -2S_n \deriv[2]{S_n}{z}+4\mu S_n \deriv{S_n}{z}. \label{eq3.23}
 \end{gather}
 Then we substitute \eqref{SnSn4} into \eqref{eq3.23} and obtain, after simplification,
 \begin{align*}
&[2n(n+1)-4\sigma]z^{2\sigl} g^2 +\cdots\\
 &\quad= -8\sigl z^{2\sigl-1}g \gz+8\mu \sigl z^{2\sigl-1} g^2-(8\sigl+2) z^{2\sigl} g \gzz +8\sigl z^{2\sigl} \left(\gz\right)^{2} +4\mu z^{2\sigl}g \gz+\cdots.
 \end{align*}
 Comparing coefficients of $z^{2\sigl-1}$ in the resulting polynomials, we obtain
\[
 8\sigl a_0 a_1-8\mu \sigl a_0^2=0.
\]
 This implies part (a). Next we compare coefficients of $z^{2\sigl}$ to get
\[
 {[2n(n+1)-4\sigl] a_0^2 = (16\sigl+4)\mu a_0 a_1-(32\sigl+4)a_0 a_2.}
\]
Since $a_1=\mu a_0$, we deduce that
\[
 a_2 =\frac{\mu^2(8\sigl+2)-n(n+1)+2\sigl}{2(8\sigl+1)} a_0.
\]
 Now since $n=\ell +|\mu|$ and $2\sigl=\ell (\ell +1)$,
\[
 n(n+1)=\mu^2+(2\ell +1)|\mu|+\ell (\ell +1),
\]
 while $8\sigl+1=(2\ell +1)^2$.
 Therefore, part (b) is valid.
\end{proof}

 \begin{Theorem} \label{thm4.6}
 Let $\mu\in \bfZ $. Then for all $n>|\mu|$, with ${n\geq 1}$:
 \begin{enumerate}\itemsep=0pt
 \item[$(a)$] for $S_n(z;\mu)$, $z=0$ is a root of order $\frac12(n-|\mu|)(n-|\mu|+1)$;
 \item[$(b)$] $S_n(z;\mu)$ is a monic polynomial of degree $\frac12n(n+1)$;
 \item[$(c)$] all other roots of $S_{n}(z;\mu)$ are simple.
 \end{enumerate}
 \end{Theorem} 

\begin{figure}[t!]
\begin{tabular}{ccc}
\fig{48mm}{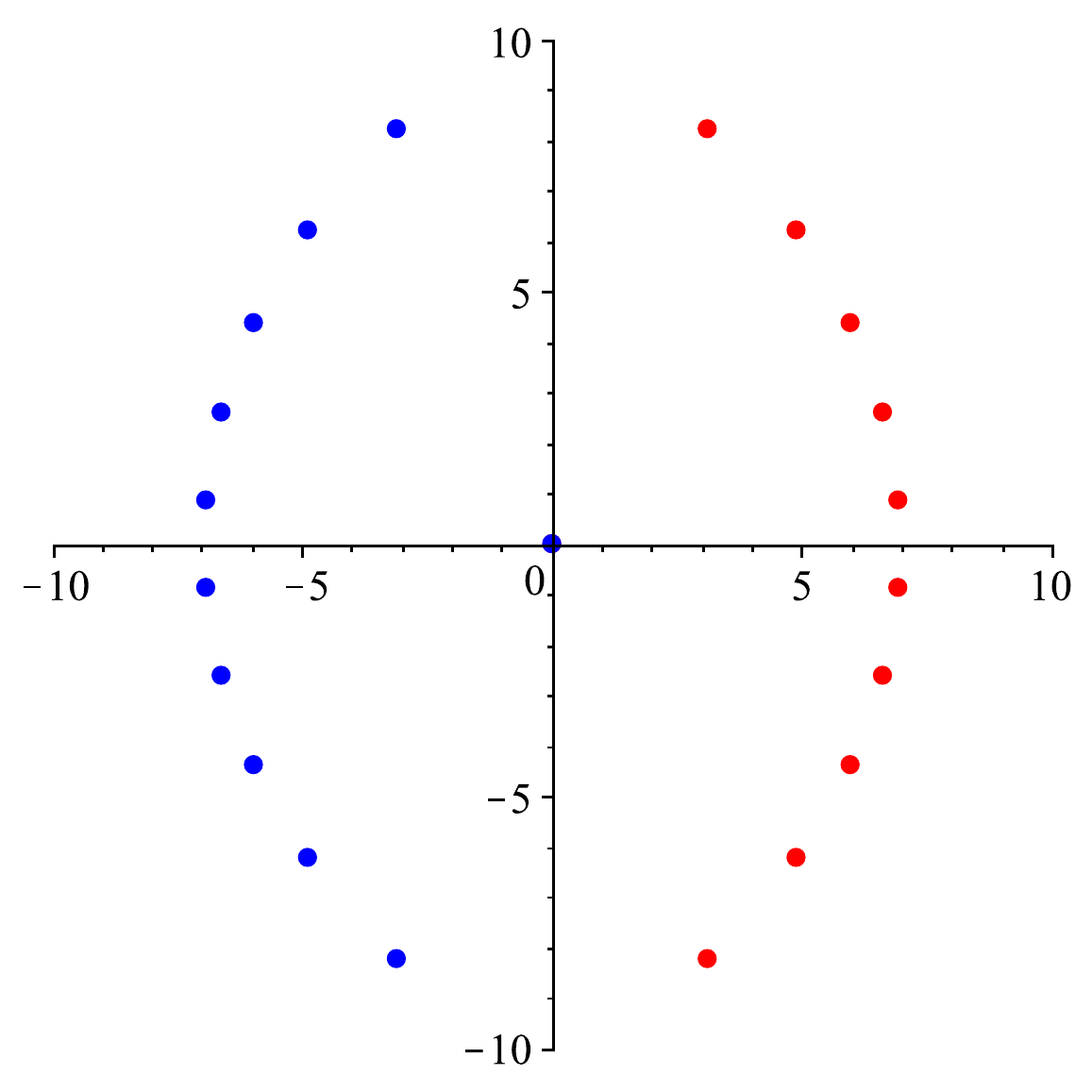} & \fig{48mm}{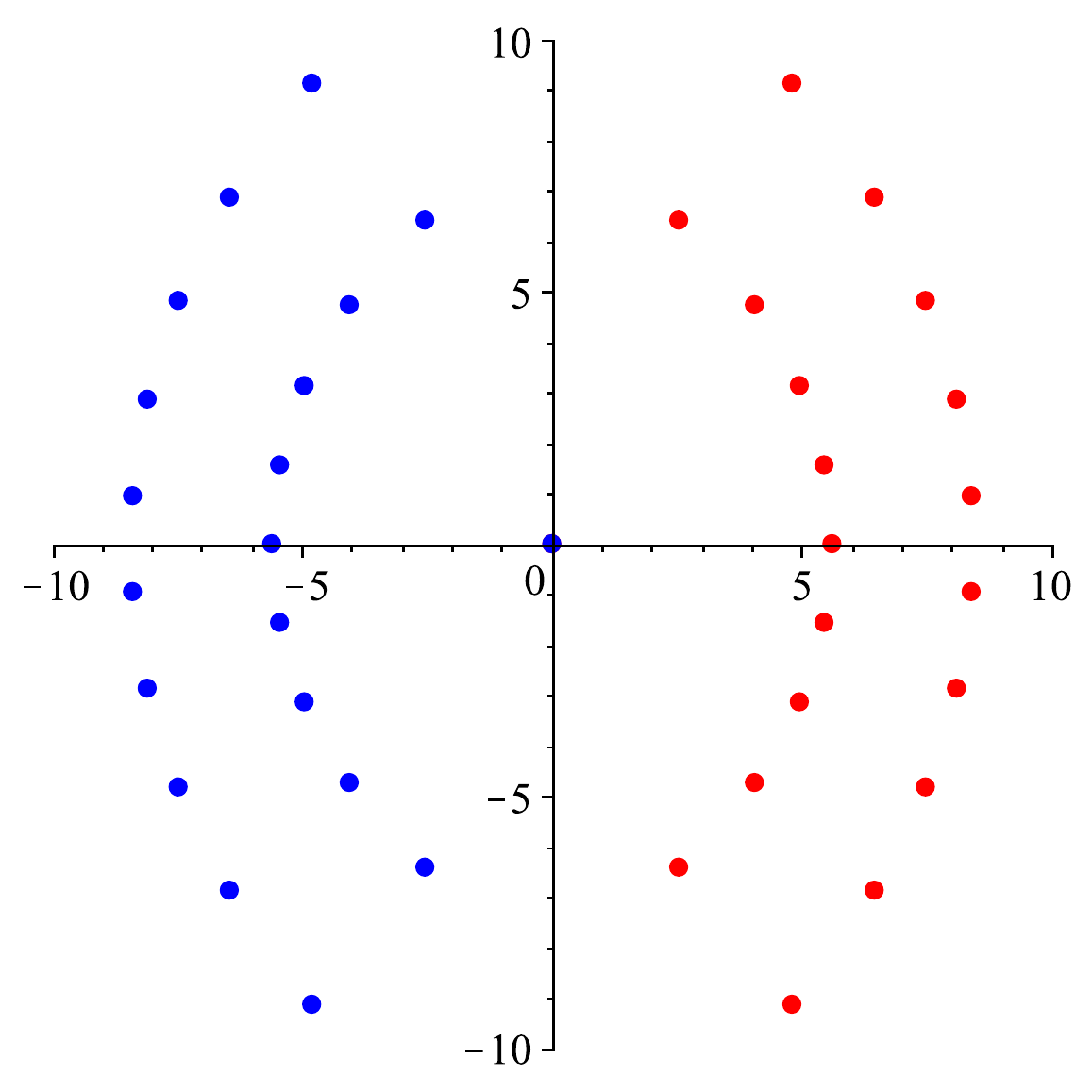} & \fig{48mm}{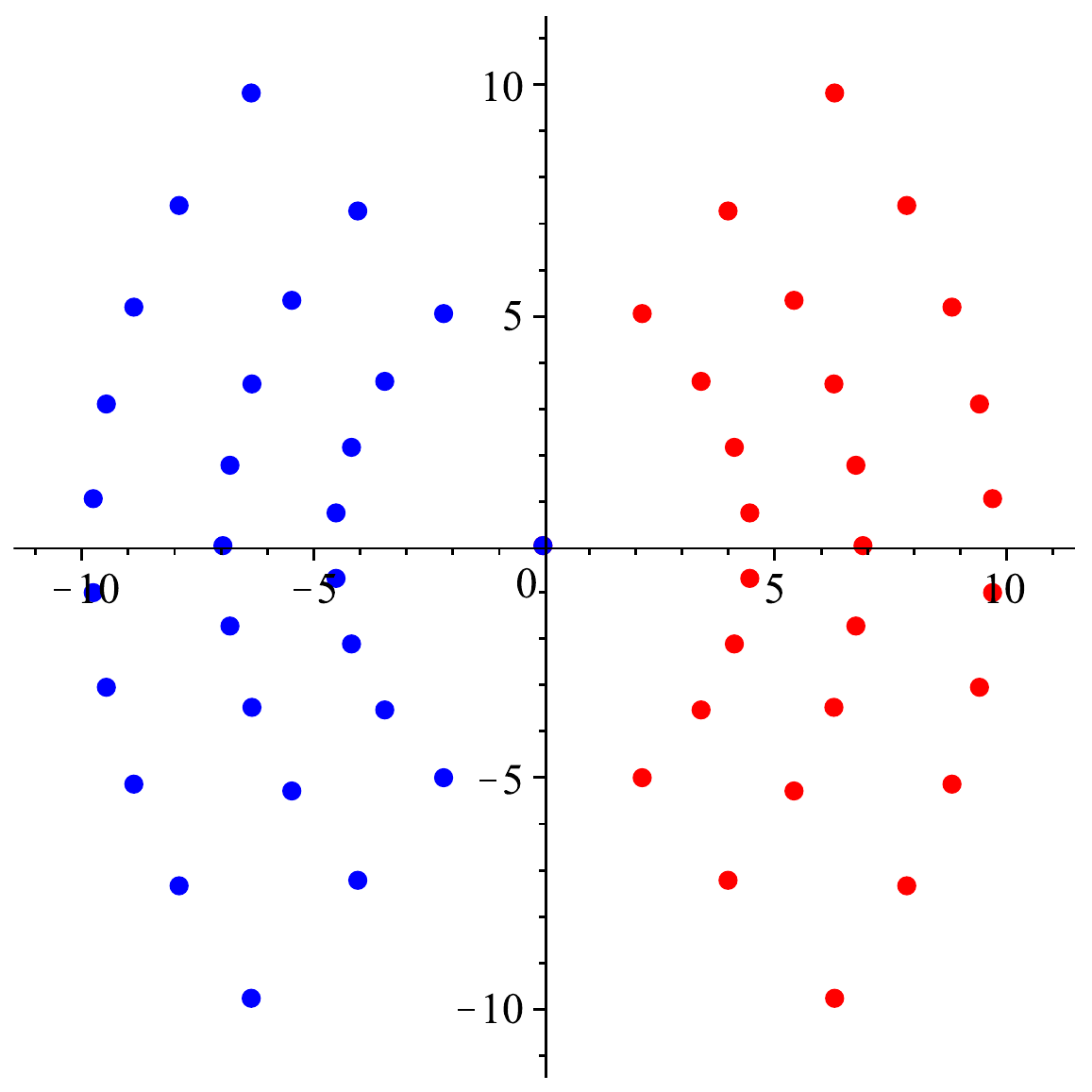}\\
 $\mu=\pm1$ & $\mu=\pm2$ & $\mu=\pm3$ \\[10pt]
\fig{48mm}{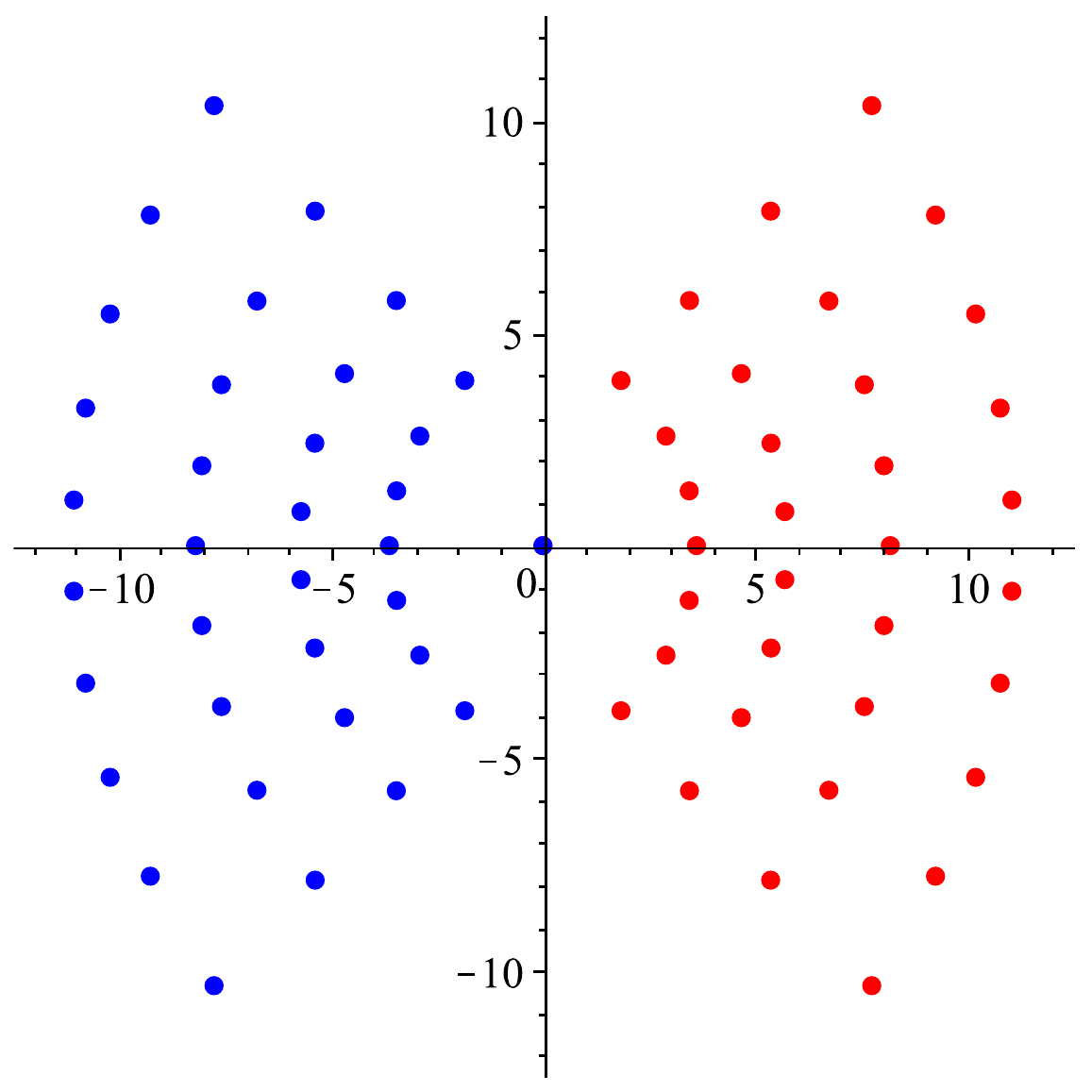} & \fig{48mm}{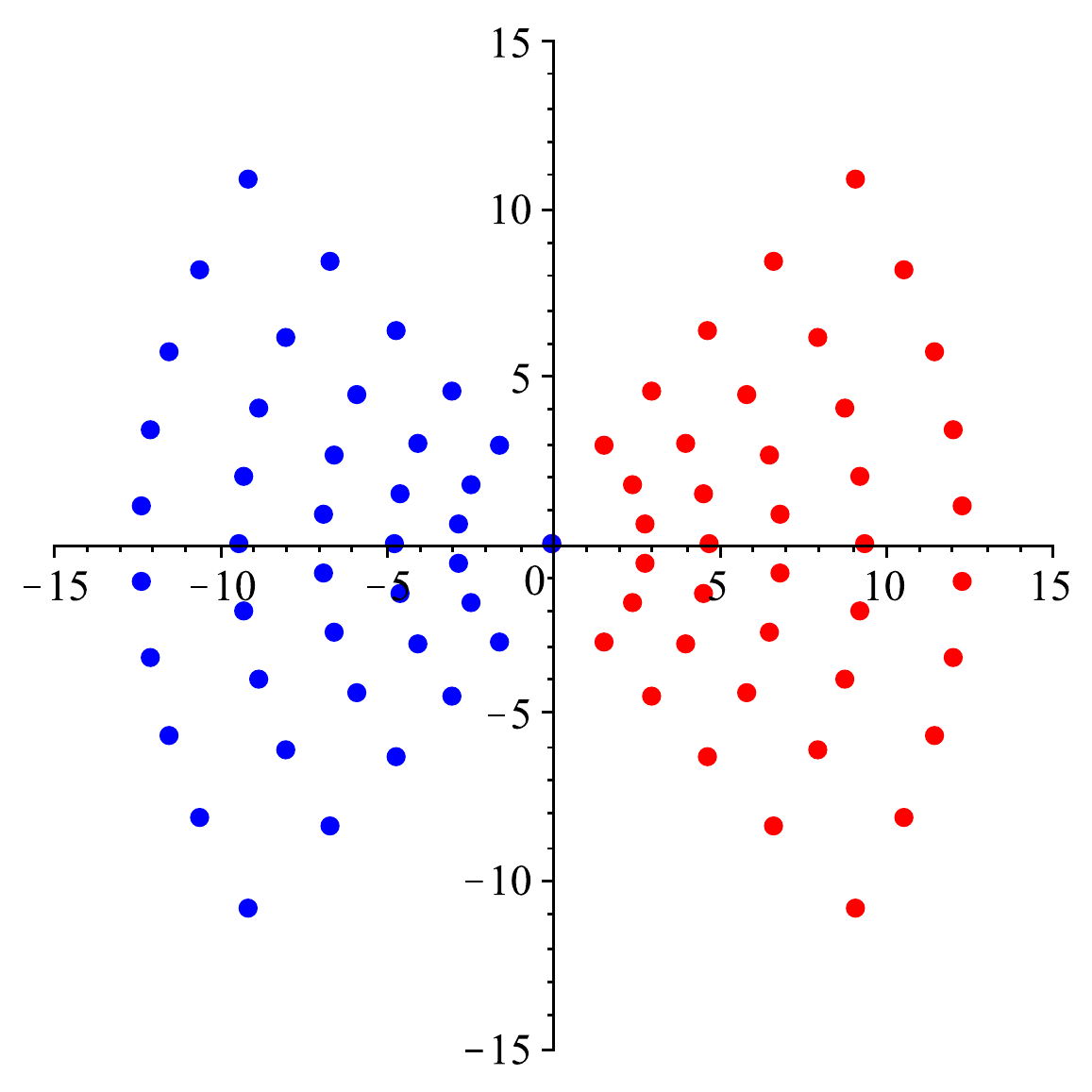} &
\fig{48mm}{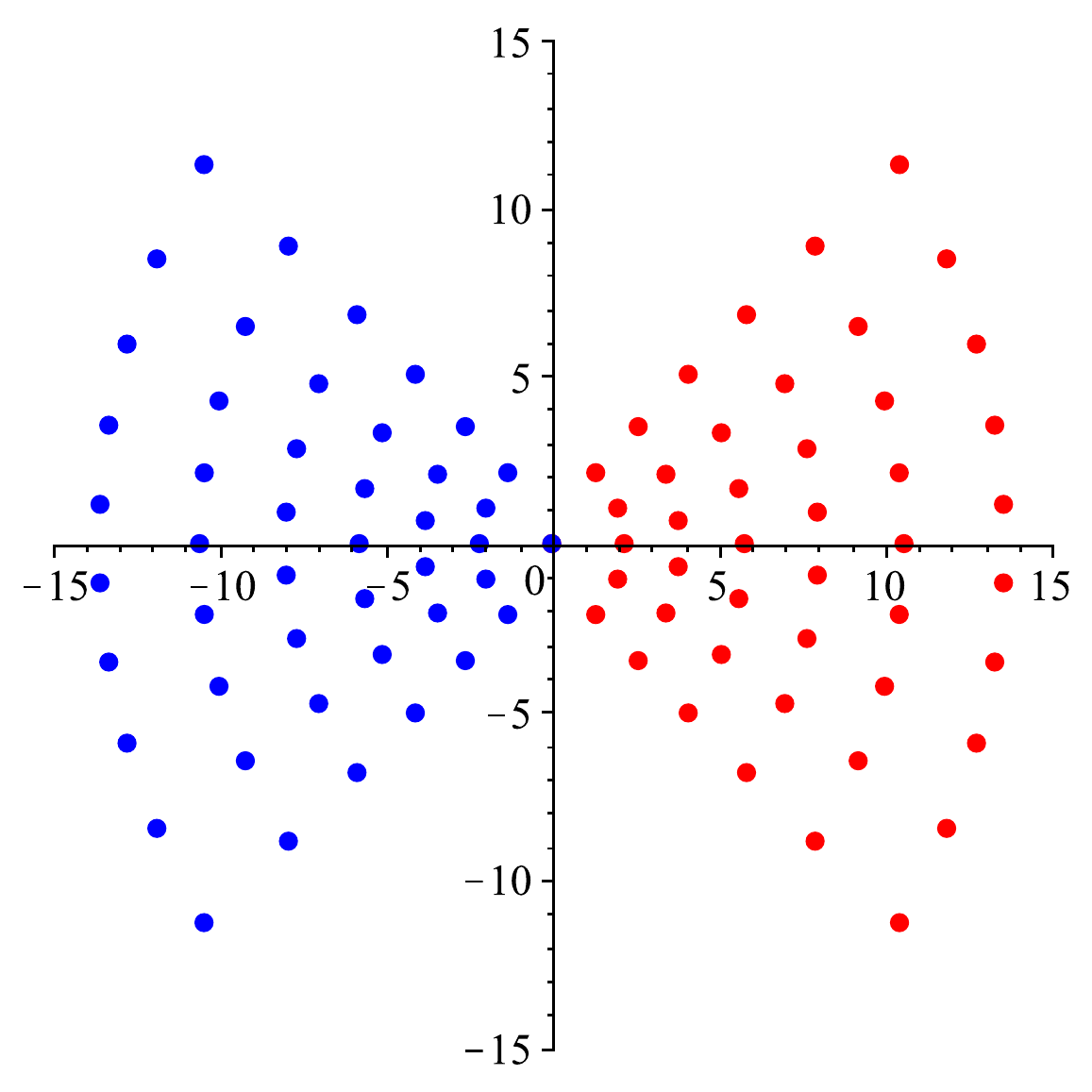}\\
 $\mu=\pm4$ & $\mu=\pm5$ & $\mu=\pm6$\\[10pt]
\fig{48mm}{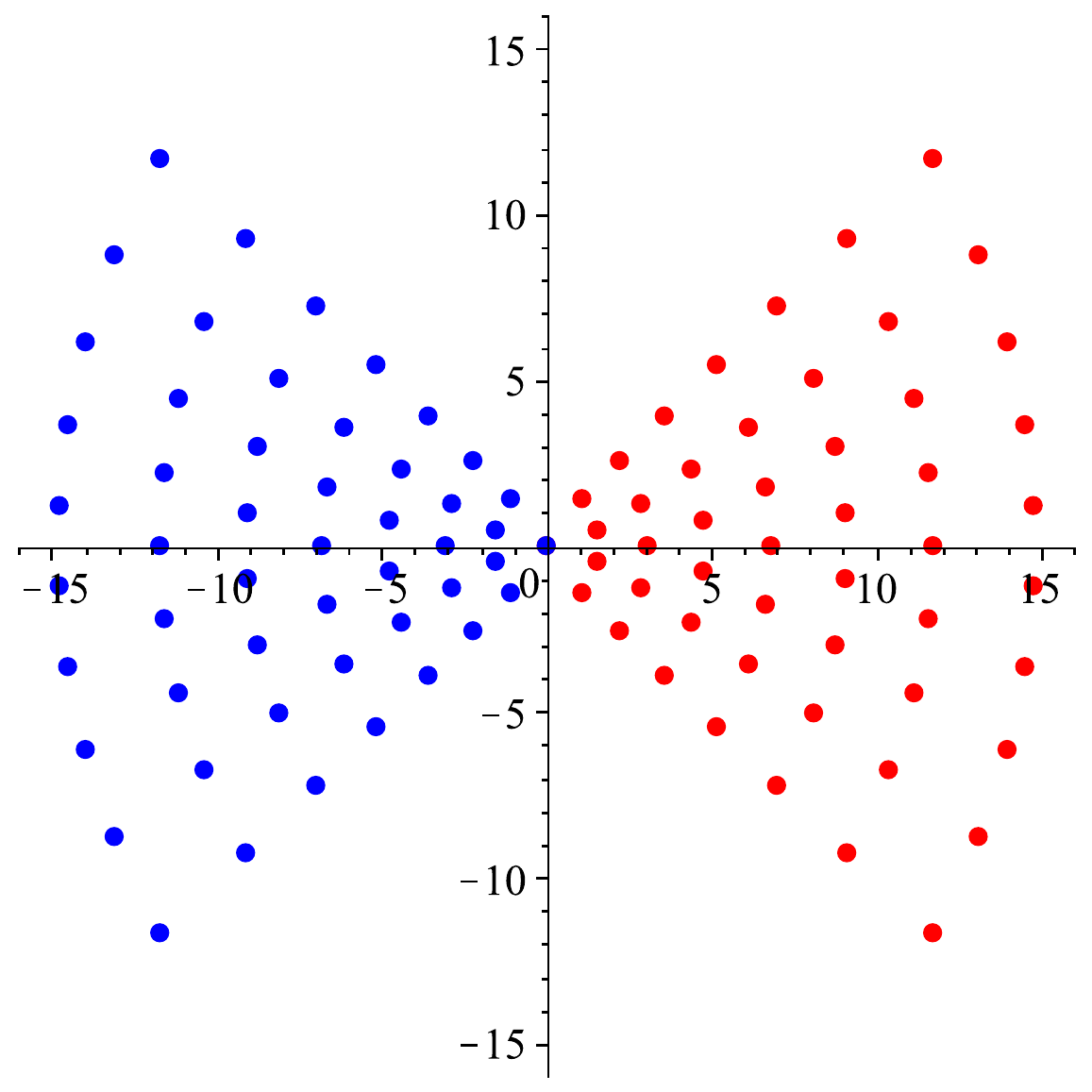} & \fig{48mm}{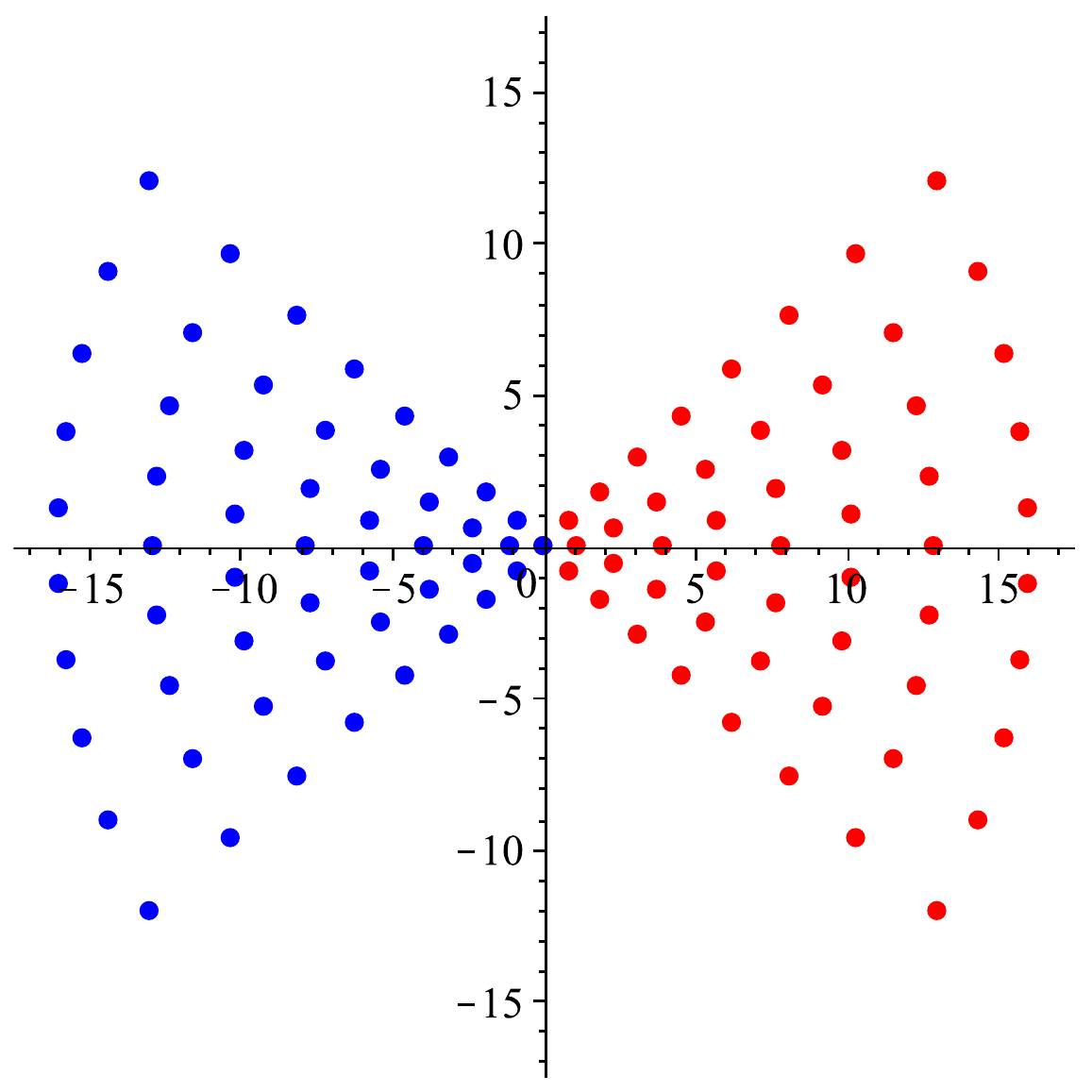} &
\fig{48mm}{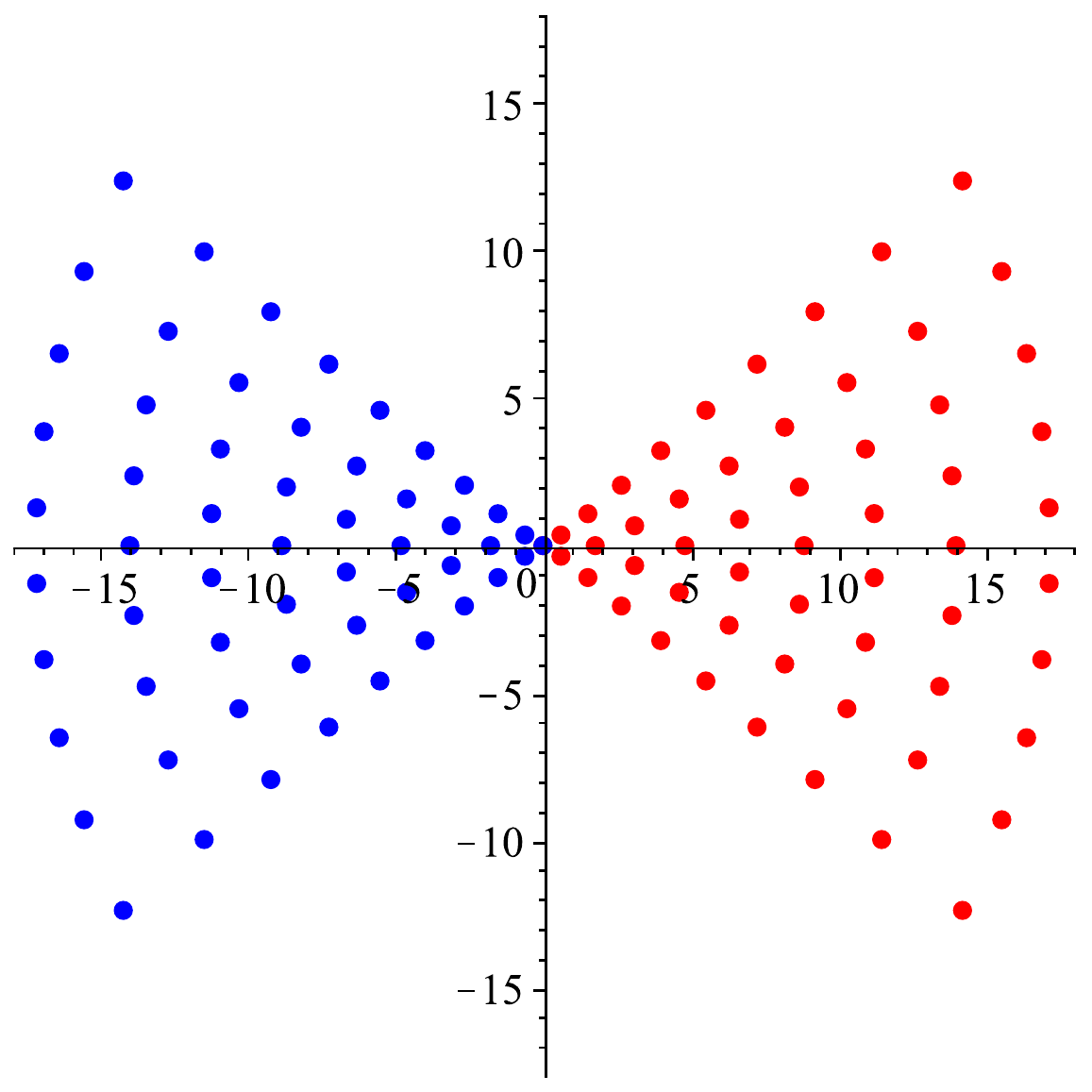}\\
 $\mu=\pm7$ & $\mu=\pm8$ & $\mu=\pm9$
 \end{tabular}
 \caption{Plots of $S_{10}(z;\mu)$ with $\mu=m$ (blue) and $\mu=-m$ (red), for $m=1,2,\dots,9$.} \label{fig:Thm46a}
 \end{figure}

\begin{figure}[t]
 \begin{tabular}{ccc}
\fig{48mm}{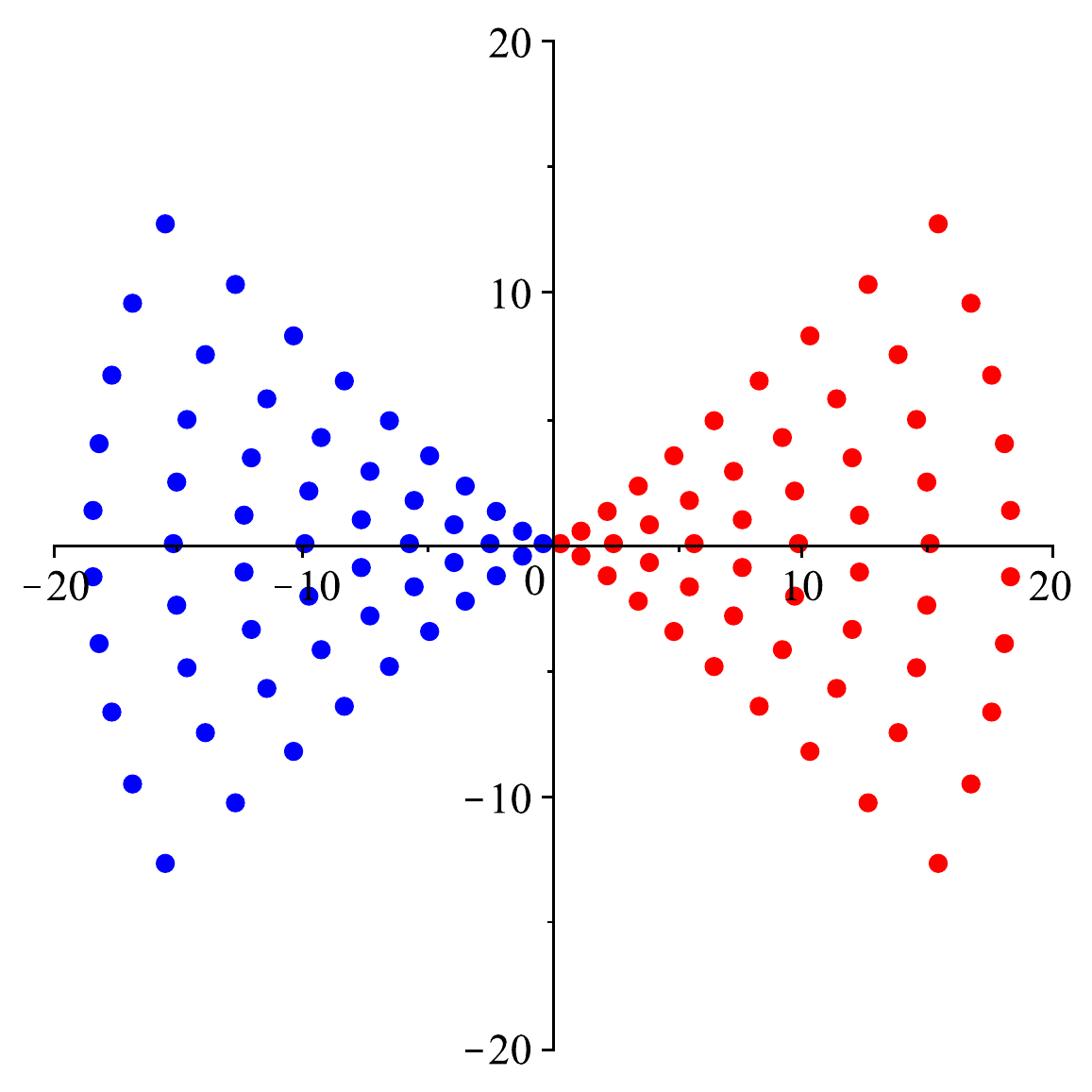} & \fig{48mm}{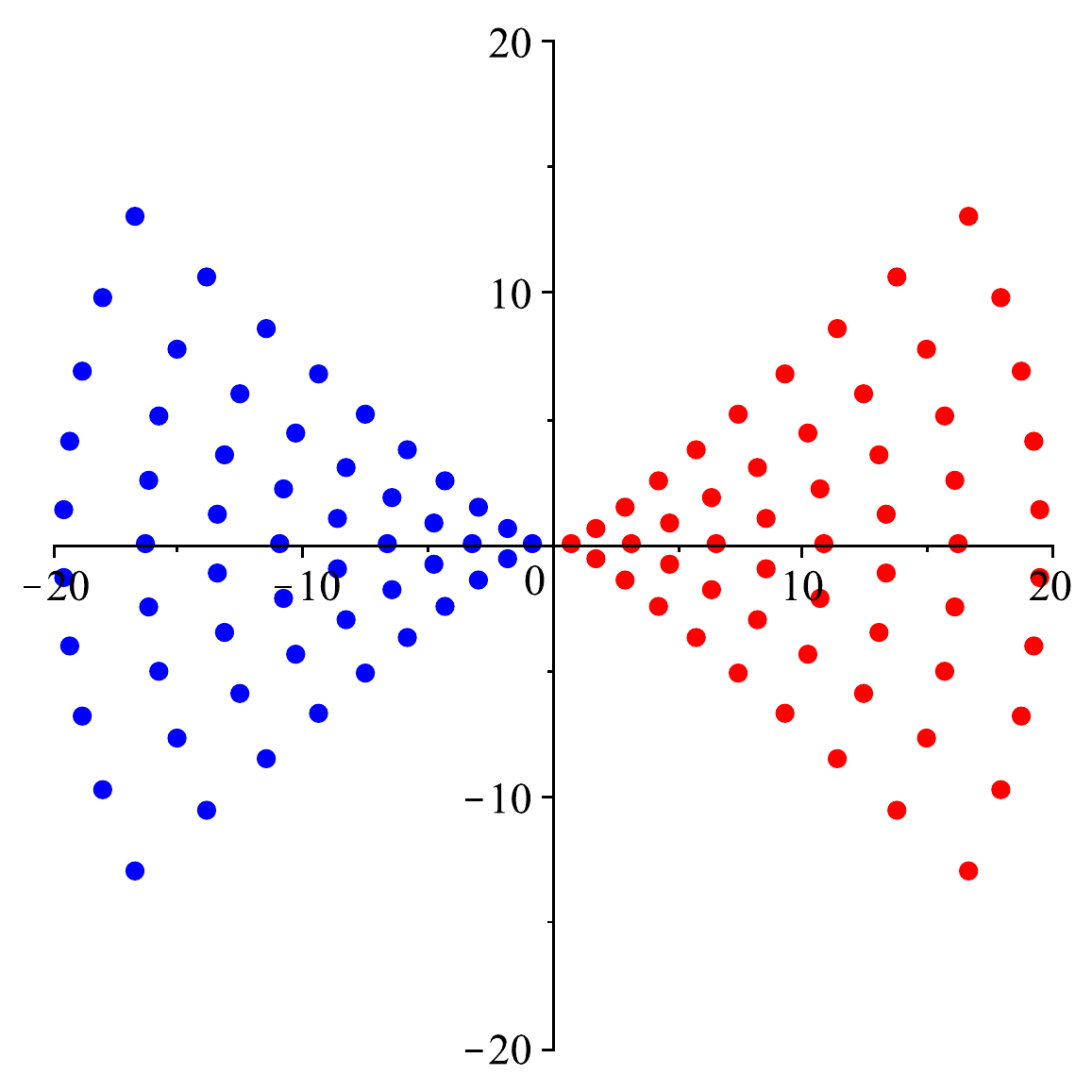} & \fig{48mm}{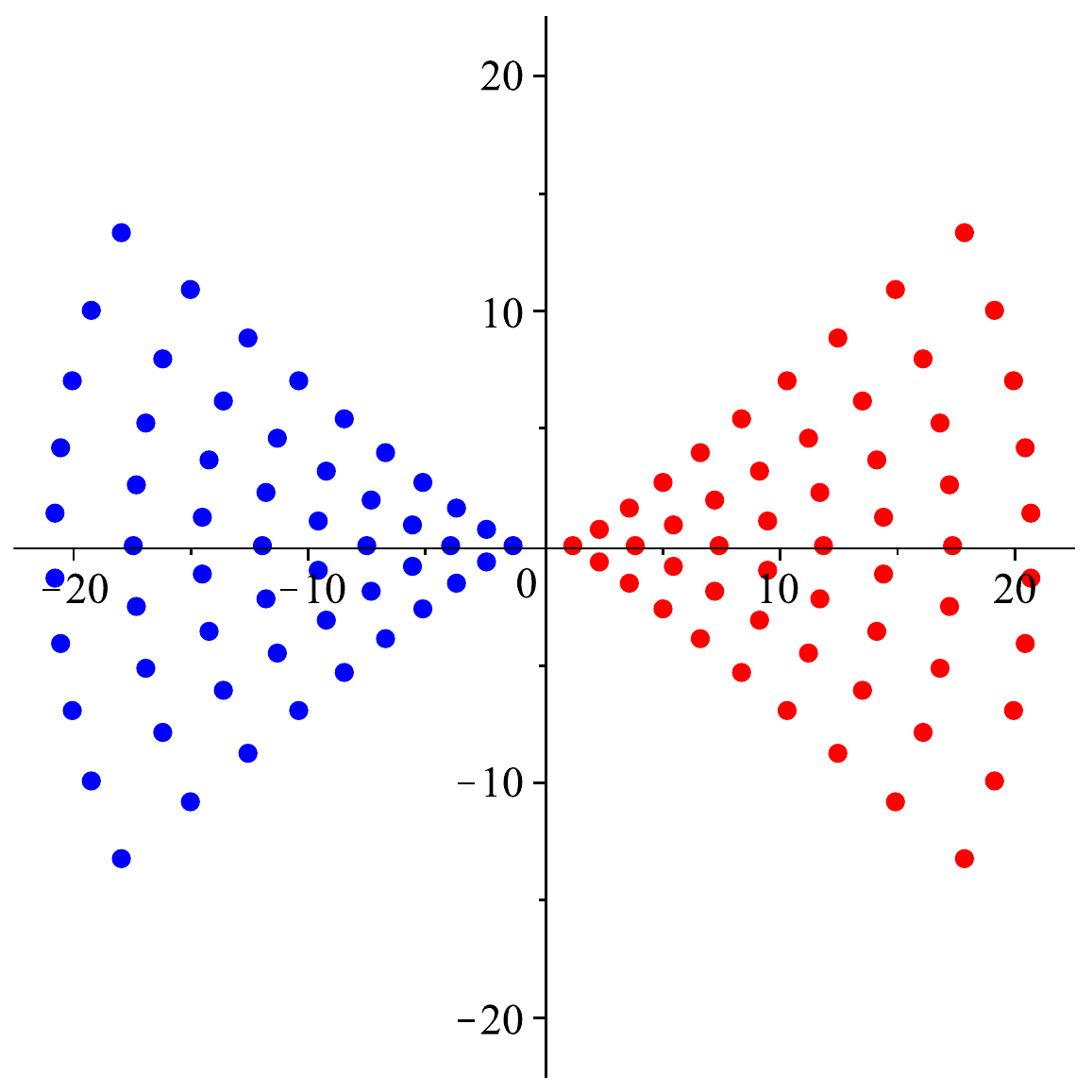}\\
 $\mu=\pm10$ & $\mu=\pm11$ & $\mu=\pm12$ \\[10pt]
\fig{48mm}{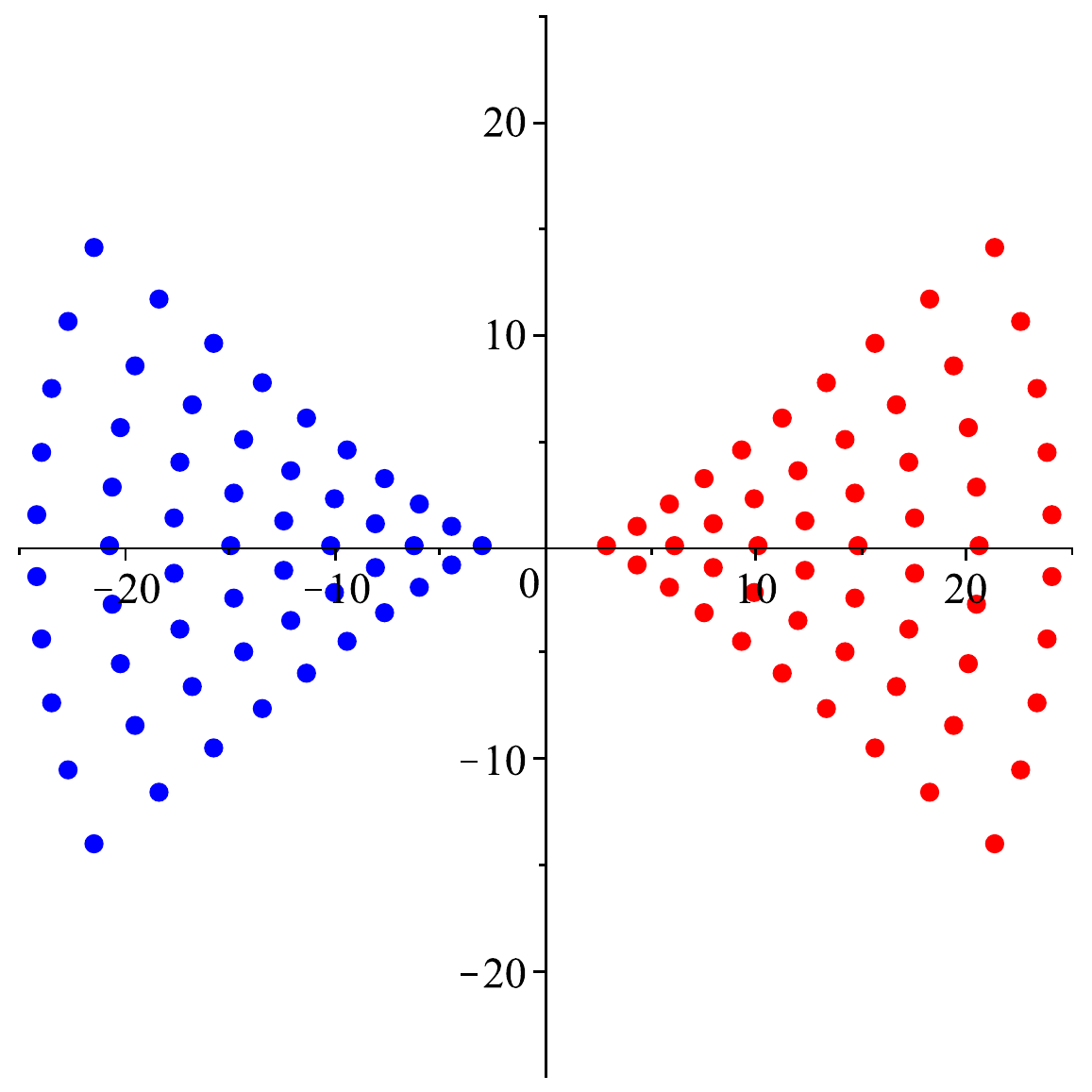} & \fig{48mm}{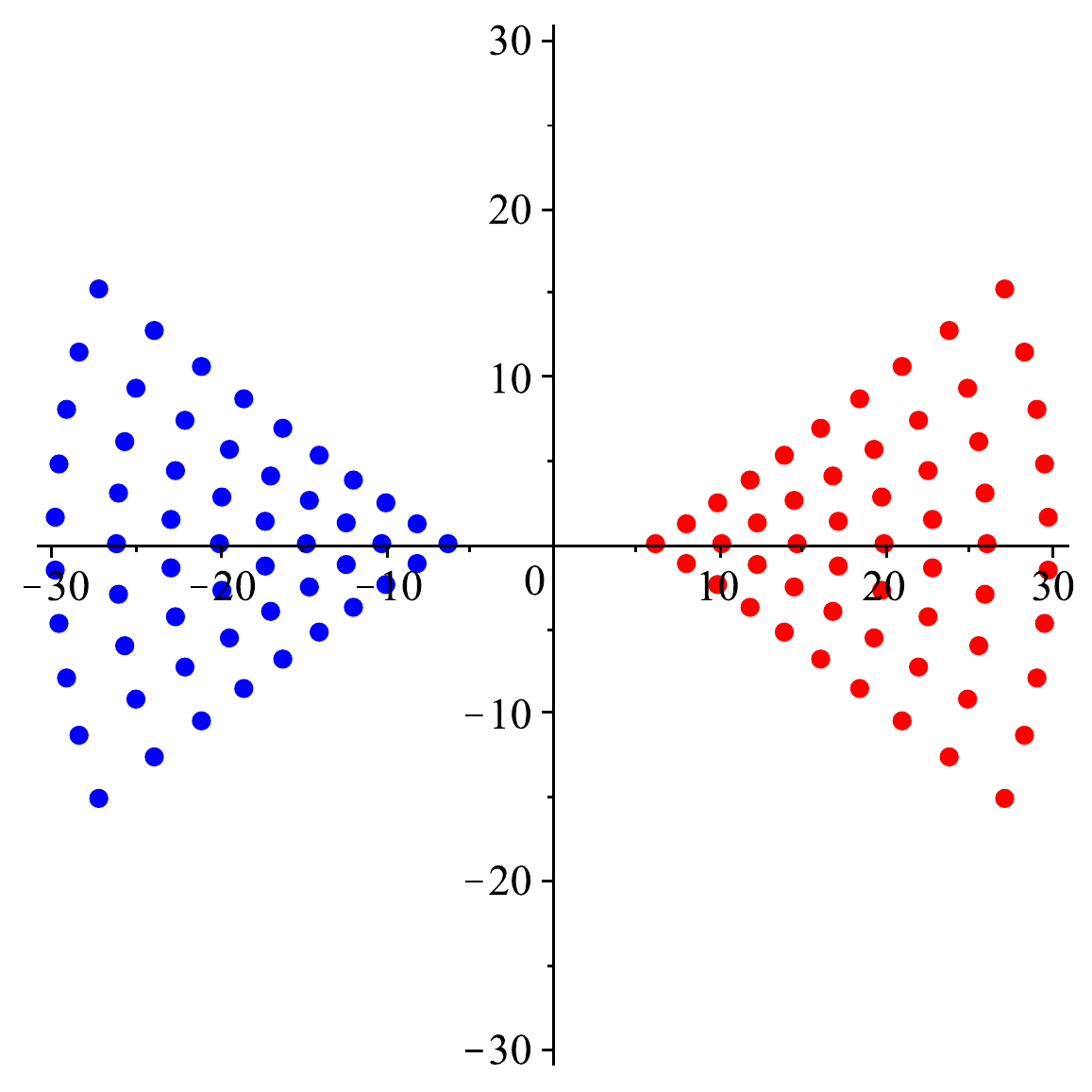} &
\fig{48mm}{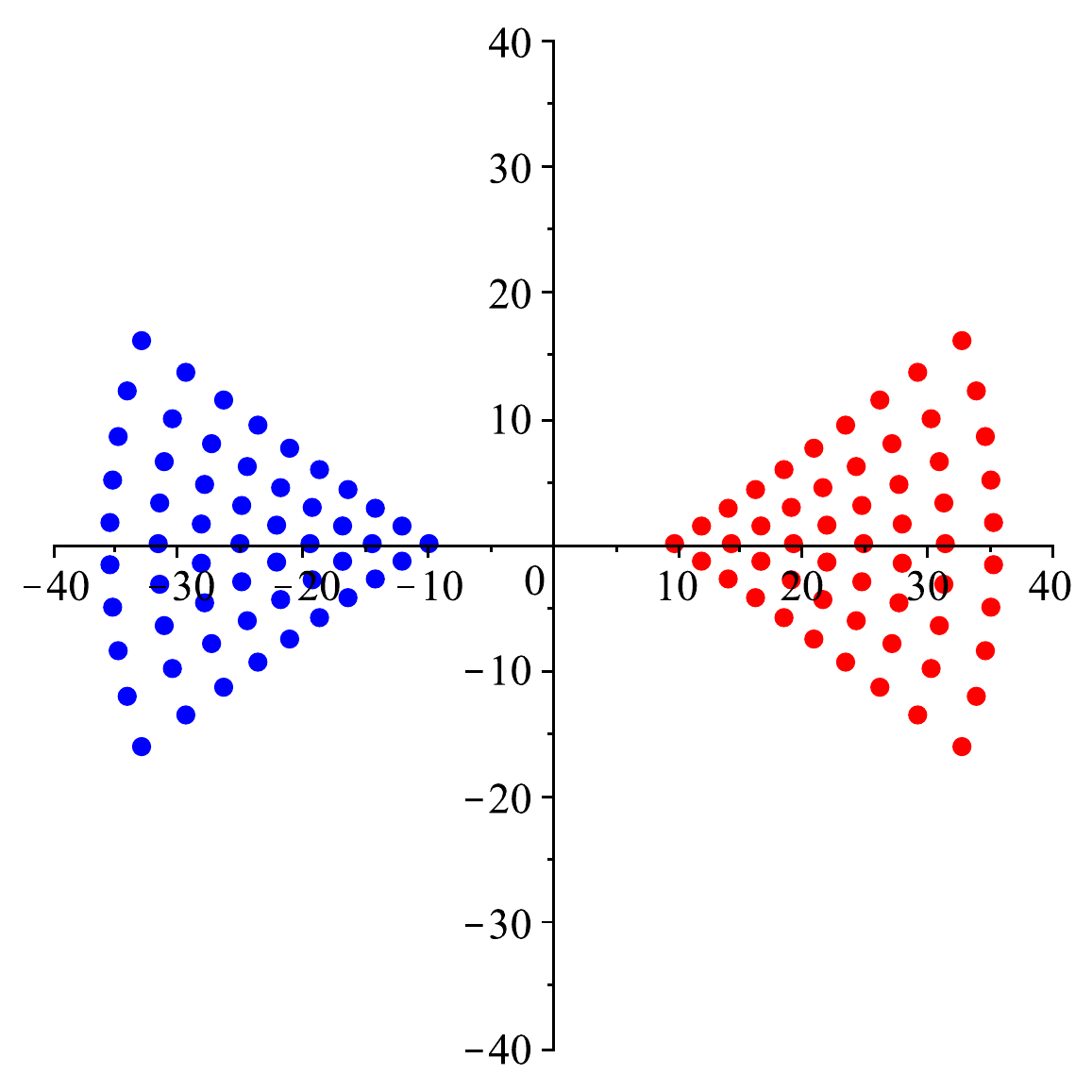}\\
 $\mu=\pm15$ & $\mu=\pm20$ & $\mu=\pm25$
 \end{tabular}
 \caption{Plots of $S_{10}(z;\mu)$ with $\mu=m$ (blue) and $\mu=-m$ (red), for $m=10, 11, 12, 15, 20, 25$. These show that for $|\mu|\geq n$, the roots of $S_{n}(z;\mu)$ have a ``triangular structure'' and lie in the region $\operatorname{Re}(z)<0$ for $\mu>0$ and $\operatorname{Re}(z)>0$ for $\mu<0$. Further as $\mu$ increases, the triangular regions move away from the imaginary axis.} \label{fig:Thm46b}
 \end{figure}

 \begin{Remark}\samepage\quad
 \begin{enumerate}\itemsep=0pt
 \item[(1)] Thus when $n$ is large, $S_n(z;\mu)$ has $\mathcal{O}\big(n^2\big)$ roots, counted according to multiplicity. But if~$\mu\in \bfZ$, then most roots are located
 at $z=0$, while there are only $\mathcal{O}(n)$ non-zero roots, and all of them are simple roots. This {explains} the phenomenon that when $\mu\in\bfZ$, the roots
 and poles of the rational solution $w_n$ are unusually fewer than the other $\mu$'s nearby, as observed in \cite{refBM,refBMS}.
\item[(2)] Theorem \ref{thm4.6} is illustrated in Figure \ref{fig:Thm46a}, where plots of $S_{10}(z;\mu)$ with $\mu=m$ (blue) and ${\mu=-m}$ (red), for $m=1,2,\dots,9$. {Contrast} this to Figure \ref{fig:Thm46b}, where plots of $S_{10}(z;\mu)$ with $\mu=m$ (blue) and $\mu=-m$ (red), for $m=10, 11, 12, 15, 20, 25$. These show that for $|\mu|\geq n$, the roots of $S_{n}(z;\mu)$ have a ``triangular structure'' and lie in the region $\operatorname{Re}(z)<0$ for $\mu>0$ and $\operatorname{Re}(z)>0$ for $\mu<0$. Further as $\mu$ increases, the triangular regions move away from the imaginary axis.
\end{enumerate}
 \end{Remark}

 \begin{proof}
 Part (c) follows from the proof of Theorem \ref{thm4.4}. For parts (a) and (b), the case when~${\mu=0}$ is simple, recall \eqref{Sn0}.
 In general, {fix any $\mu\in \bfZ\setminus\{ 0\}$ and let $m=|\mu|$.} By~Corollary~\ref{th3.2},
 $z=0$ is a root of $S_{m+1}(z;\mu)$.
 Observe that by {\eqref{eq3.3}, \eqref{eq3.4}} and \eqref{eq3.13},
\[
 \deriv{S_{m+1}}{z}(0;\mu) =\phi_{m+1}'=\frac{\phi_m\left[\phi_m+2\mu \phi_m'-2\phi_m''\right]}{\phi_{m-1}}=\frac{\phi_m^2}{\phi_{m-1}}(2m+1)\neq 0,
\]
 because $\phi_m''=\big(m^2-m\big)\phi_m$. Thus $z=0$ is a simple root of $S_{m+1}(z;\mu)$, and we may write $S_{m+1}(z;\mu)=z g_{1}(z)$, which by Theorem~\ref{th2.12} is a polynomial.
 Let $ {g_{1}(z)=\sum_{j=0}^k a_j^{1} z^j}$ be a~polynomial of degree $k= {\frac12(m+1)(m+2)}-1$, with nonzero roots.

 Now we apply the induction hypothesis on $n\geq m+1$. Let $\ell=n-m$, and
\[
 S_{n-1}(z;\mu)=z^{\sigma_{\ell-1}}g_{\ell-1}(z),\qquad S_{n}(z;\mu)=z^{\sigma_{\ell}}g_{\ell}(z),
\]
 where $g_{\ell-1}(z;\mu)$ and $g_{\ell}(z;\mu)$ are polynomials with nonzero roots {and $\sigma_{\ell}=\frac12\ell(\ell+1)$}.
 Then by~\eqref{eq:Srecrel},
\[
 {z^{\sigma_{\ell-1}} g_{\ell-1}S_{n+1}(z;\mu)} = z^{2\sigma_{\ell}}\left(\mu g_{\ell}^2-g_{\ell}\deriv{g_{\ell}}{z}\right)
 +z^{2\sigma_{\ell}+1}\bigg[\left(\deriv{g_{\ell}}{z}\right)^{2}-g_{\ell}\deriv[2]{g_{\ell}}{z}+g_{\ell}^2\bigg].
\]
 Let $a_j^{\ell}$ be the coefficients of $g_{\ell}$. By Lemma \ref{thm4.5}, $a_1^{\ell}=\mu a_0^{\ell}$. So we may write
 $S_{n+1}(z;\mu)=z^{\sigma_{\ell+1}}g_{\ell+1}$, where
 \begin{align*}
 a_0^{\ell+1} a_0^{\ell-1} = 2\mu a_0^{\ell}a_1^{\ell}-4 a_0^{\ell}a_2^{\ell}+\big(a_0^{\ell}\big)^{2}
 = a_0^{\ell}\big(2\mu a_1^{\ell}+a_0^{\ell}-4 a_2^{\ell}\big)
 = \big(a_0^{\ell}\big)^{2} \left(1+\frac{2m}{2\ell +1}\right).
 \end{align*}
 So {$a_0^{\ell+1}=g_{\ell+1}(0)$} is nonzero, and {the function $g_{\ell+1}(z)$, which is a rational function at its initial appearance}, does not have $z=0$ as a root.

 Next we show that $g_{\ell+1}(z)$ is a polynomial. From the proof of Theorems~\ref{thm4.4} and~\ref{th2.12}\,(b),
 we know that all nonzero roots of $S_{n}(z;\mu)$ and $S_{n-1}(z;\mu)$ are simple and not common. Furthermore, we still have $S_{n-1} \,|\, {\left[-z \L_z(S_n)+(z+\mu)S_n^2\right]}$, where
\[
 \Delta:=-z \L_z(S_n)+(z+\mu)S_n^2=z\bigg[\left(\deriv{S_n}{z}\right)^{2}-S_n\deriv[2]{S_n}{z}\bigg]-S_n \deriv{S_n}{z}+(z+\mu)S_n^2.
\]
We conclude that $g_{\ell-1}$ divides $\Delta$, which implies that $g_{\ell+1}$ is indeed a polynomial. It means $S_{n+1}=z^{\sigma_{\ell+1}}g_{\ell+1}(z)$ is indeed a polynomial.
Consequently, parts (a) and (b) follow by induction.\looseness=-1
\end{proof}

\section{Conclusions}\label{sec5}
{We have given a direct algebraic proof that the nonlinear recurrence relation \eqref{eq:Srecrel} generates polynomials $S_{n}(z;\mu)$, rather than rational functions without direct resort to the $\tau$-function theory of
\p\ equations. However we critically needed a higher order equation derived from the corresponding $\sigma$-equation, which seems to be inevitable in the nonlinear scenario.
We believe that the method can be developed to apply to the fifth \p\ equation ($\PV$) as well, though we shall not pursue this further here.}

\appendix
\section{About the coalescence limit}
\begin{Lemma}
\label{tha.1}
The sequence of functions
\[
 R_n(\zeta,\ep):= \ep^{n(n+1)/2} S_n\left(\frac{\zeta}{\ep}+\frac{4}{\ep^3},-\frac{4}{\ep^3}\right)
\]
 are all polynomials in $\ep$.
 \end{Lemma}
\begin{proof}
From \eqref{eq:Srecrel}, we write
 \begin{equation}
 S_{n+1}S_{n-1}=-(z+\mu)\big(S_nS_n''-(S_n')^2\big)-S_nS_n'+(z+\mu)S_n^2+\mu\big(S_nS_n''-(S_n')^2\big)
 \label{eqa.1}
 \end{equation}
 with $S_0=S_{-1}=1$. It is easy to see from Theorems \ref{th2.12} and~\ref{thm4.6} that each $S_n$ is a~polynomial in $\zeta=z+\mu$, as well as a polynomial in $\mu$.
 Furthermore,
\[
 \deg (S_n,\zeta)=\frac12 n(n+1)=\deg (S_n,\mu).
\]
 Now let
\[
 V_n\big(\zeta,\ep^{-1}\big):=S_n\left(\frac{\zeta}{\ep}+\frac{4}{\ep^3}, -\frac{4}{\ep^3}\right).
\]
 We claim that $V_n$ is a polynomial in $\ep^{-1}$, and $\deg \big(V_n,\ep^{-1}\big)=\frac12 n(n+1)$,
 so that each $R_n$ defined above, as a rational function, is indeed a polynomial in $\ep$.

Rewrite \eqref{eqa.1} as
 \begin{align}
 V_{n-1}V_{n+1} ={}& -\frac{\zeta}{\ep}\bigg\{V_n \deriv[2]{V_n}{\zeta} - \left(\deriv{V_n}{\zeta}\right)^{2}\bigg\}-V_n \deriv{V_n}{\zeta}+ \frac{\zeta}{\ep}V_n^2\nonumber
 \\
 &-\frac{4}{\ep^3}\bigg\{V_n \deriv[2]{V_n}{\zeta} - \left(\deriv{V_n}{\zeta}\right)^{2}\bigg\}.
 \label{eqa.2}
 \end{align}
 Hence $V_0=V_{-1}=1$, and
 \[V_1=\frac{\zeta}{\ep},\qquad V_2=\frac{\zeta^3+4}{\ep^3}, \]
 and so on. It is trivial to show that $V_n$ is a polynomial in $\zeta$, with $\deg(V_n,\zeta)=\frac12 n(n+1)$. Moreover, inductively, the right-hand side of \eqref{eqa.2} is a polynomial in $\ep^{-1}$, with degree $n^2+n+1$, and the leading coefficient of $\ep^{-(n^2+n+1)}$ involves $\zeta^{n^2+n+1}$, and is so nonzero. By induction hypothesis, $\deg\big(V_{n-1},\ep^{-1}\big)=\frac12 n(n-1)$, and $V_{n-1}$ divides the expression on the right-hand side. Therefore, we have
 $V_{n+1}$ is also a polynomial in $\ep^{-1}$, and
\[
 \deg\big(V_{n+1},\ep^{-1}\big)=\frac12 (n+1)(n+2).
\]

We emphasize that in the above argument, the terms in $S_n(z,\mu)$ can achieve maximum power of $\ep^{-1}$ only at those terms involving $\zeta^{n(n+1)/2-3k}$, so that for each derivative with respect to $\zeta$, the power of $\ep^{-1}$ will decrease by one. Also the coefficient at the maximum power of $\ep^{-1}$ does not vanish because of the expression $\zeta^{n^2+n+1}$ in the third term above.

 The proof is now complete.
 \end{proof}

\subsection*{Acknowledgements}
The authors are deeply indebted to the anonymous referees for their very careful and constructive review work.
PAC thanks Thomas Bothner, Alfredo Dea\~no, Clare Dunning, Marco Fasondini, Kerstin Jordaan, Ana Loureiro and Walter Van Assche for their helpful comments and illuminating discussions and also the Department of Applied Mathematics, National Sun Yat-sen University, Kaohsiung, Taiwan and the Department of Mathematics, National Taiwan University, Taipei, Taiwan, for their hospitality during his visit where some of this work was done.
CKL also thanks Peter Miller, Yik-Man Chiang and Guofu Yu for stimulating discussions.
CKL is partially supported by National Science and Technology Council (formerly Ministry of Science and Technology), Taiwan.

\pdfbookmark[1]{References}{ref}
\LastPageEnding

\end{document}